\documentclass[a4paper,12pt]{article}
\usepackage{pdfsync}
\usepackage{amsfonts}
\usepackage{amssymb}
\usepackage{amsmath}
\usepackage{mathptmx}
\usepackage{cancel}
\usepackage{pstricks}
\usepackage{psfrag}
\usepackage{mathrsfs}
\usepackage{graphicx}
\usepackage{pgf,tikz}
\usepackage{algorithm}
\usepackage{algpseudocode}
\usepackage{soul}
\def\be{\begin{equation}}
\def\ee{\end{equation}}
\def\bea{\begin{eqnarray}}
\def\eea{\end{eqnarray}}
\def\beann{\begin{eqnarray*}}
\def\eeann{\end{eqnarray*}}

\newcommand{\rank}{{\rm rank}}

\def\ns{\hspace{-1mm}}
\newcommand{\real}{{\mathbb{R}}}

\def\spacingset#1{\def\baselinestretch{#1}\small\normalsize}
\spacingset{1.30}

\newtheorem{lemma}{Lemma}
\newtheorem{theorem}{Theorem}
\newtheorem{remark}{Remark}
\newtheorem{corollary}{Corollary}
\newtheorem{proposition}{Proposition}

\newtheorem{example}{Example}[section]

\def\be{\begin{equation}}
\def\ee{\end{equation}}
\def\bea{\begin{eqnarray}}
\def\eea{\end{eqnarray}}
\def\beann{\begin{eqnarray*}}
\def\eeann{\end{eqnarray*}}

\def\ns{\hspace{-1mm}}

\def\proof{\noindent{\bf{\em Proof:}\ \ }}
\def\QED{\mbox{\rule[0pt]{1.5ex}{1.5ex}}}
\def\endproof{\hspace*{\fill}~\QED\par\endtrivlist\unskip}

\def\second{{\prime \prime}}

\newcommand{\ima}{\operatorname{im}}

\newcommand{\defi}{\stackrel{\text{\tiny def}}{=}}

\definecolor{Royalblue}{cmyk}{1,0.30,0.2,0.2}

\newcommand{\complex}{{\mathbb{C}}}

\def\gU{{\cal U}}
\def\gV{{\cal V}}

\def\gX{{\cal X}}

\def\bmat{\left[ \begin{array}}
\def\emat{\end{array} \right]}

\def\bmat{\left[ \begin{array}}
\def\emat{\end{array} \right]}
\def\bsmat{\left[ \begin{smallmatrix}}
\def\esmat{\end{smallmatrix} \right]}

\def\gU{{\cal U}}

\def\gV{{\cal V}}

\def\gX{{\cal X}}

\def\sss{\scriptscriptstyle}

\newcommand{\spanR}{\operatorname{span}}
\begin{document}
\begin{titlepage}
\title{\vspace{-5mm}
The Discrete-Time Generalized 
Algebraic Riccati Equation: Order Reduction and Solutions' Structure\thanks{Partially supported by the Australian Research Council under the grant FT120100604 and by the
Italian Ministry for Education and Research (MIUR) under PRIN grant n. 20085FFJ2Z.  }\vspace{10mm}}
\author{Lorenzo Ntogramatzidis$^\star$ \and Augusto Ferrante$^\dagger$}
\date{$^\star${\small Department of Mathematics and Statistics,\\[-2pt]
         Curtin University, Perth (WA), Australia \\[-2pt]
         {\tt L.Ntogramatzidis@curtin.edu.au}\\[8pt]
        $^\dagger$Dipartimento di Ingegneria dell'Informazione,\\[-2pt]
         Universit\`a di Padova, via Gradenigo, 6/B -- 35131 Padova, Italy\\[-2pt]
         {\tt augusto@dei.unipd.it} }
        }%
\thispagestyle{empty} \maketitle \thispagestyle{empty}
\begin{abstract}%
In this paper we discuss how to decompose the constrained generalized discrete-time algebraic Riccati equation arising in optimal control and optimal filtering problems {into two parts corresponding to an additive decomposition $X=X_0+\Delta$ of each solution $X$: The first part is an explicit expression of the addend $X_0$ which is common to all solutions, and does not depend on the particular $X$. The second part can be 
 either a reduced-order 
discrete-time  {\em regular} algebraic Riccati equation whose  associated closed-loop matrix is non-singular, or a symmetric Stein equation. }\end{abstract}

\begin{center}
\begin{minipage}{14.2cm}
\vspace{2mm}
\end{minipage}
\end{center}
\thispagestyle{empty}
\end{titlepage}
%

\section{Introduction}
\label{secintro}

This paper is concerned with the {following relations}
\bea
X  \ns&\ns = \ns&\ns  A^{\sss \top} X A  -  (A^{\sss \top} X B  +  S)(R  +  B^{\sss \top} X B)^\dagger(S^{\sss \top} +B^{\sss \top} X A)  +  Q, \label{gdare} \\
 &   &  \qquad \ker (R+B^{\sss \top}\,X\,B) \subseteq \ker (A^{\sss \top}\,X\,B+S) \label{kercond}
\eea
where the symbol $\dagger$ denotes the Moore-Penrose pseudo-inverse operation.\footnote{We recall that given an {\em arbitrary} matrix $M\in \real^{h \times k}$, there exists a unique matrix $M^\dagger\in \real^{k \times h}$  that satisfies the following four properties: {\bf (1)} $M\,M^\dagger\,M=M$; {\bf (2)} $M^\dagger\,M\,M^\dagger=M^\dagger$; {\bf (3)} $M^\dagger\,M$ is symmetric; {\bf (4)} $M\,M^\dagger$ is symmetric.
 By definition, the matrix $M^\dagger$ is  the {\em Moore-Penrose pseudo-inverse}  of the matrix  $M$.} Equation (\ref{gdare}) subject to the constraint (\ref{kercond}) arises for example in {discrete-time} LQ problems -- see \cite{Rappaport-S-71} and \cite{Ferrante-N-12} for the finite and infinite-horizon cases, respectively. Here, $A  \in \real^{n \times n}$, $B \in \real^{n \times m}$, $Q\in \real^{n \times n}$, $S \in \real^{n \times m}$ and $R \in \real^{m \times m}$ are such that the {\em Popov matrix} $\Pi$ satisfies
\be
\label{popov}
\Pi \defi \bmat{cc} Q & S \\[-1mm] S^{\sss \top} & R \emat =\Pi^{\sss \top} \ge 0.
\ee
The set of matrices $\Sigma=(A,B;\Pi)$ is often referred to as {\em Popov triple}, and  (\ref{gdare}) is known  as the {\em generalized discrete-time algebraic Riccati equation} GDARE($\Sigma$). This equation, together with the additional constraint (\ref{kercond}), is usually referred to as {\em constrained generalized discrete-time algebraic Riccati equation}, and it is herein denoted by CGDARE($\Sigma$). This
equation generalizes the standard {\em discrete-time algebraic Riccati equation} DARE($\Sigma$) 
\bea
X  \ns&\ns = \ns&\ns  A^{\sss \top} X A  -  (A^{\sss \top} X B  +  S)(R  +  B^{\sss \top} X B)^{-1}(S^{\sss \top} +B^{\sss \top} X A)  +  Q, \label{dare} 
\eea
as the natural equation arising in LQ optimal control and filtering problems.
In fact, it is only when the underlying linear  system | obtained by a full-rank factorization  $\Pi=\bsmat C^{\sss \top} \\[1mm] D^{\sss \top} \esmat [\begin{array}{cc} C & D \end{array} ]$ and considering a system described by the quadruple $(A,B,C,D)$ | is left invertible that the standard DARE($\Sigma$) admits solutions. The dynamic optimization problem, however, may still admit solutions in the more general setting where the underlying linear system is not  left-invertible. In these cases, however, the standard DARE($\Sigma$) does not admit solutions and the correct equation that must be used to address the original optimization problem is the CGDARE($\Sigma$), see e.g. \cite{Ferrante-N-12}. As discussed in \cite[Chapt. 6]{Abou-Kandil-FIJ-03}, these general situations are particularly relevant in the context of stochastic control problems, see also \cite{Damm-H-01,Freiling-H-04} and the references cited therein.
On the other hand, whenever the standard DARE($\Sigma$) admits solutions, the set of its solutions coincides with the set of solutions of CGDARE($\Sigma$), so that the latter is a genuine generalization of the former
(here and in the rest of the paper, we are only considering {\em symmetric} solutions $X$ both for the DARE($\Sigma$) and the CGDARE($\Sigma$)). 

In the literature, several efforts have been devoted by many authors to the task of reducing the order and difficulty of the standard DARE($\Sigma$) by means of different techniques, \cite{Mita-85,Fujinaka-A-87,Fujinaka-CS-99,Hansson-H-99,Ferrante-04,Ferrante-W-07}. This interest is motivated by the fact that the standard DARE($\Sigma$) is richer than the structure of its continuous-time counterpart, the continuous-time algebraic Riccati equation. In particular, in \cite{Ferrante-04} a method was presented which, differently from earlier contributions presented on this topic, aimed at iteratively decomposing DARE($\Sigma$)  into a trivial part and a reduced DARE whose associated closed-loop matrix is non-singular. 
 The subsequent contribution \cite{Ferrante-W-07} achieves a similar goal by avoiding the need for an iterative procedure. 

The development of reduction procedures for generalized Riccati equations has received much less attention in the literature. This is in part likely to be due to the technical difficulties associated with generalized Riccati equations in the discrete time. In \cite{Ferrante-04}, a hint is given on how the iterative reduction detailed therein could be extended to the case of an equation in the form (\ref{gdare}), provided that the attention is restricted to the set of positive semidefinite solutions, {for which} condition (\ref{kercond}) is automatically satisfied. 
 On the other hand, {CGDARE($\Sigma$) may well admit solutions that are not} positive semidefinite, see e.g. \cite{Ferrante-N-12,Ferrante-N-13-1}.
In \cite{Hansson-H-99}, a Riccati equation in the form of a CGDARE($\Sigma$) is considered, and a reduction technique is proposed to the end of computing the stabilizing solution of CGDARE($\Sigma$). The main goal of this paper is to combine the generality of the framework considered in \cite{Hansson-H-99} with the ambition of achieving a reduction for the entire set of solutions of  CGDARE($\Sigma$). This task is accomplished by developing an iterative procedure that is similar in spirit to that of \cite{Ferrante-04}, but which presents a richer and more articulated structure. Indeed, not only do several technical difficulties and structural differences arise in extending the results of 
\cite{Ferrante-04} to the case of CGDARE($\Sigma$) when the set of solutions is not restricted to semidefinite ones, but also, 
 differently from the iterations needed in  \cite{Ferrante-04}, which are essentially performed via changes of coordinates in the state space, in the general case of a CGDARE($\Sigma$), it is necessary to also resort to changes of coordinates in the input space. The problem of obtaining a systematic procedure to decompose generalized Riccati equations into a trivial part and a reduced, ``well-behaved'', part described by a {\em regular} DARE (or at times, differently from the standard case, by a symmetric Stein equation), becomes much more interesting and challenging in the case of generalized Riccati equations.
 Our reduction method is based on the computation of null spaces of given matrices so that  it can be easily implemented in a software procedure that uses only standard linear algebra procedures
 which are robust and available in any numerical software package. Therefore
a relevant outcome of the presented procedure is what we believe to be the first systematic numerical procedure to compute the solutions of CGDARE.

 %
 %

\section{Problem formulation and preliminaries}\label{sec2}
First, in order to simplify the notation, for any $X=X^{\sss \top}\in  \real^{n \times n}$ we define the matrices
\beann
\begin{array}{ccccllllcccccccc}
R_{\sss X} & \defi  & R+B^{\sss \top}\,X\,B &&& G_{\sss X} & \defi  & I_{\sss m}-R_{\sss X}^\dagger R_{\sss X}\\
S_{\sss X} & \defi & A^{\sss \top}\,X\,B+S &&& 
K_{\sss X} & \defi  & R_{\sss X}^\dagger S_{\sss X}^{\sss \top} &&& 
A_{\sss X}&  \defi & A-B\,K_{\sss X}
\end{array}
\eeann
so that (\ref{kercond}) in CGDARE($\Sigma$) can be written concisely as $\ker R_{\sss X} \subseteq \ker S_{\sss X}$. The term $R_{\sss X}^\dagger R_{\sss X}$ is the orthogonal projector that projects onto $\ima R_{\sss X}^\dagger=\ima R_{\sss X}$ so that $G_{\sss X}$ is  the orthogonal projector that projects onto $\ker R_{\sss X}$. Hence,
 $\ker R_{\sss X}=\ima G_{\sss X}$.

As already mentioned, in this paper we present a procedure that reduces CGDARE($\Sigma$) to another discrete-time algebraic Riccati equation with the same structure { but smaller order and in which}  both $A_{\scriptscriptstyle 0} \defi A-B\,R^\dagger\,S^{\sss \top}$ and $R$ are non-singular. On the other hand,  this means that the Riccati equation thus obtained is indeed a standard DARE, i.e., it has the structure shown in (\ref{dare}), as the following result shows. 
\begin{proposition}
Suppose that the matrix $R$ is non-singular, and let $X=X^\top$ be any symmetric solution of CGDARE($\Sigma$). Then $R_{\sss X} = R  +  B^{\sss \top} X B$ is non-singular.\\[-3mm]
\end{proposition}
\proof As shown in \cite[Lemma 4.1]{Ferrante-N-12}, for any symmetric solution $X=X^\top$ of CGDARE($\Sigma$) the inclusion $\ker R_{\scriptscriptstyle X} \subseteq \ker R$ holds. As a consequence, if $R$ is non-singular, its null-space $\ker R$ is zero, and therefore so is the null-space of $R_{\scriptscriptstyle X}$. This is equivalent to the fact that $R_{\scriptscriptstyle X}$ is non-singular.
\endproof

%

The reduction technique presented in this paper can also be viewed from the perspective of the so-called
extended symplectic pencil $N_{\sss \Sigma}-z\,M_{\sss \Sigma}$, where
\beann
M_{\sss \Sigma}\defi\left[ \begin{array}{ccc} 
I_{\sss n} & 0 & 0 \\
0 & -A^{\sss \top} & 0 \\ 
0 & -B^{\sss \top} & 0 \end{array} \right] \qquad \textrm{and} \qquad
N_{\sss \Sigma}\defi \left[ \begin{array}{ccc} 
A & 0 & B \\
Q & -I_{\sss n} & S \\ 
S^{\sss \top} & 0 & R \end{array} \right].  
\eeann
The case in which the matrix pencil $N_{\sss \Sigma}-z\,M_{\sss \Sigma}$ is regular (i.e., if there exists $z \in \complex$ such that $\det (N_{\sss \Sigma}-z\,M_{\sss \Sigma}) \neq 0$) corresponds to the case in which CGDARE($\Sigma$) is indeed a DARE($\Sigma$), whereas the one in which $N_{\sss \Sigma}-z\,M_{\sss \Sigma}$ is singular (i.e., the determinant of $N_{\sss \Sigma}-z\,M_{\sss \Sigma}$ is the zero polynomial) corresponds to a case in which DARE($\Sigma$) does not admit solutions.
%
It is shown in \cite{Ferrante-04} for DARE($\Sigma$) and in \cite{Ferrante-N-13-2} for CGDARE($\Sigma$) that if $A_{\sss X}$ is singular, the Jordan structure of $A_{\sss X}$ associated with the eigenvalue $\lambda=0$ is completely determined by 
$N_{\sss \Sigma}-z\,M_{\sss \Sigma}$, and is independent of the particular solution $X$ of DARE($\Sigma$) or CGDARE($\Sigma$). It is shown in \cite{Ferrante-04} that in the case where the matrix pencil $N_{\sss \Sigma}-z\,M_{\sss \Sigma}$ is regular   | or, equivalently,  the CGDARE($\Sigma$) and  the standard DARE($\Sigma$) have the same solutions| the following statements are equivalent:
\begin{description}
\item{\bf (1)} $N_{\sss \Sigma}$ is singular;
\item{\bf (2)} $N_{\sss \Sigma}-z\,M_{\sss \Sigma}$ has a generalized eigenvalue at zero;
\item{\bf (3)}  there exists a solution $X$ of CGDARE($\Sigma$) such that the corresponding closed-loop matrix $A_{\sss X}$ is singular;
\item{\bf (3$^\prime$)}   for any solution $X$ of CGDARE($\Sigma$), the corresponding closed-loop matrix $A_{\sss X}$ is singular;
\item{\bf (4)} at least one of the two matrices $R$ and $A_{\sss 0}=A-B\,R^\dagger\,S^{\sss \top}$ is singular.
\end{description}
The case where the matrix pencil $N_{\sss \Sigma}-z\,M_{\sss \Sigma}$ is possibly singular was investigated in \cite{Ferrante-N-13-2}, where it was proved that in this more general case these four facts are not equivalent. In particular, {\bf (1)} is not equivalent to {\bf (2)}. Moreover, in the case where $N_{\sss \Sigma}-z\,M_{\sss \Sigma}$ is singular, {\bf (1)} and {\bf (3)} are not equivalent, nor are {\bf (3)} and {\bf (4)}. 
However, 
it was shown in \cite[Lemma 3.1]{Ferrante-N-13-2}
that it is still true that  {\bf (1)} is equivalent to {\bf (4)}.
Furthermore,  
 it is shown in \cite[Proposition 3.4]{Ferrante-N-13-2} that $r\defi \rank R_{\sss X}$ is constant  for any solution $X$ of 
CGDARE($\Sigma$), and that
 $A_{\sss X}$ is singular if and only if at least one of the following two conditions holds:
 (i) $\rank R < r=\rank R_{\sss X}$ and (ii) $A_{\sss 0}=A-B\,R^\dagger\,S^{\sss \top}$ is singular. It is clear that this condition reduces to {\bf (4)} in the case where $R_{\sss X}$ is invertible, i.e., in the case where $N_{\sss \Sigma}-z\,M_{\sss \Sigma}$ is regular.
Notice also that since both conditions are independent of the particular solution $X$ of the CGDARE($\Sigma$), the singularity of the closed-loop matrix $A_X$ is invariant with respect to the particular solution $X$.
 
  To summarize, in the case where the matrix pencil $N_{\sss \Sigma}-z\,M_{\sss \Sigma}$ is singular, the following statements are equivalent:
\begin{description}
\item{\bf (1$^\prime$)} $N_{\sss \Sigma}$ is singular;
\item{\bf (2$^\prime$)} at least one of the two matrices $R$ and $A_{\sss 0}=A-B\,R^\dagger\,S^{\sss \top}$ is singular;
\end{description}
and the following statements are equivalent:
\begin{description}
\item{\bf (1$^\second$)}   there exists a solution $X$ of CGDARE($\Sigma$) such that  the corresponding closed-loop matrix   $A_{\sss X}$ is singular;
\item{\bf (2$^\second$)}   for any solution $X$ of CGDARE($\Sigma$), the corresponding closed-loop matrix $A_{\sss X}$ is singular;
\item{\bf (3$^\second$)} at least one of the two conditions 
\begin{description}
\item{\bf (a)} $\rank R < r=\rank R_{\sss X}$; or 
\item{\bf (b)} $A_{\sss 0}=A-B\,R^\dagger\,S^{\sss \top}$ is singular;
\end{description}
is satisfied.
\end{description}
%
%
%
We recall again that in \cite[Lemma 4.1]{Ferrante-N-12} it was shown that for any solution X of CGDARE($\Sigma$), we have $\ker R_{\sss X} \subseteq \ker R$. This means that 
if $R$ is non-singular, such is also $R_{\sss X}$, and therefore the condition $\rank R < \rank R_{\sss X}$ is not satisfied. Thus, in this case,  the closed-loop matrix $A_{\sss X}$ is non-singular for some solution $X$ of the CGDARE($\Sigma$) if and only if it is non-singular for each solution $X$ of the CGDARE($\Sigma$) and this is in turn equivalent to $A_{\sss 0}$ being non-singular.


\section{Mathematical preliminaries}
We begin this section by recalling a standard linear algebra result that is used in the derivations throughout the paper.
  {
  \begin{lemma}
  \label{lem1}
  Consider $P=\left[ \begin{smallmatrix} P_{\sss 11} & P_{\sss 12} \\[1mm] P_{\sss 12}^{\sss \top} & P_{\sss 22} \end{smallmatrix} \right]=P^{\sss \top} \ge0$. Then,
    \begin{description}
  \item{\bf (i)} $\,\ker P_{\sss 12} \supseteq \ker P_{\sss 22}$; 
  \item{\bf (ii)} $\,P_{\sss 12}\,P_{\sss 22}^\dagger\,P_{\sss 22}= P_{\sss 12}$;  
  \item{\bf (iii)} $\,P_{\sss 12}\,(I-P_{\sss 22}^\dagger P_{\sss 22})=0$; 
  \item{\bf (iv)} $\,P_{\sss 11}-P_{\sss 12} P_{\sss 22}^\dagger P_{\sss 12}^{\sss \top} \ge 0$.
  \end{description}
    \end{lemma}

\ \\[-3mm]
We now generalize a well-known result of the classic Riccati theory | which essentially shows how to eliminate the cross-penalty matrix $S$ | to the case of a constrained generalized Riccati equation.

\begin{lemma}
\label{lem21}
Let $A_{\scriptscriptstyle 0}  \defi   A-B\,R^\dagger\,S^{\sss \top}$ and $Q_{\scriptscriptstyle 0}  \defi   Q-S\,R^\dagger\,S^{\sss \top}$.
Moreover, let $\Pi_{\sss 0}\defi \bsmat Q_{\sss 0} && 0 \\[1mm] 0 && R \esmat$
and $\Sigma_{\sss 0} \defi (A_{\sss 0},B,\Pi_{\sss 0})$. Then, the following statements hold true:
\begin{description}
\item{\bf (i)} CGDARE($\Sigma$) has the same set of solutions as CGDARE($\Sigma_{\sss 0}$)
\bea
X \ns&\ns=\ns&\ns   A_{\scriptscriptstyle 0}^{\sss \top} X A_{\scriptscriptstyle 0}  -  A_{\scriptscriptstyle 0}^{\sss \top} X B (R  +  B^{\sss \top} X B)^\dagger B^{\sss \top} X A_{\scriptscriptstyle 0} +  Q_{\scriptscriptstyle 0},\label{gdare0}  \\
 &   &  \qquad \ker (R+B^{\sss \top}\,X\,B) \subseteq \ker (A_{\scriptscriptstyle 0}^{\sss \top}\,X\,B);
 \label{kercond0}
\eea
\item{\bf (ii)} for any symmetric solution $X$ of CGDARE($\Sigma$), we have
{\beann
A_{\scriptscriptstyle X} \ns&\ns = \ns&\ns A_{\scriptscriptstyle 0X}   \defi  A_{\scriptscriptstyle 0}-B\, (R  +  B^{\sss \top} X B)^\dagger B^{\sss \top} X A_{\scriptscriptstyle 0};
\eeann
}
\item{\bf (iii)} $Q_{\scriptscriptstyle 0} \ge 0$.
\end{description}
\end{lemma}
\proof
We start proving {\bf (i)}. 
Inserting the expressions for $A_{\scriptscriptstyle 0}$ and $Q_{\scriptscriptstyle 0}$ into (\ref{gdare0}) yields
\bea
X \ns&\ns =\ns&\ns  
 %
 A^{\sss \top} \,X\,A
 -A^{\sss \top} \,X\,B\,R^\dagger S^{\sss \top} -S\,R^\dagger B^{\sss \top}\,X\,A
 +S\,R^\dagger B^{\sss \top}\,X\,B\,R^\dagger S^{\sss \top} \nonumber \\
\ns&\ns\ns&\ns  -A^{\sss \top} X\,B\,R_{\scriptscriptstyle X}^\dagger B^{\sss \top} X\, A
 +A^{\sss \top}  X\,B\,R_{\scriptscriptstyle X}^\dagger B^{\sss \top} X\,B\,R^\dagger S^{\sss \top} 
 +S\,R^\dagger B^{\sss \top} X\,B\,R_{\scriptscriptstyle X}^\dagger B^{\sss \top} X\,A \nonumber
 \\
 \ns&\ns\ns&\ns -S\,R^\dagger B^{\sss \top} X\,B\,R_{\scriptscriptstyle X}^\dagger B^{\sss \top} X\,B\,R^\dagger S^{\sss \top}+Q-S\,R^\dagger S^{\sss \top} \nonumber \\
 \ns&\ns =\ns&\ns A^{\sss \top} \,X\,A
 -A^{\sss \top} \,X\,B\,R^\dagger S^{\sss \top} -S\,R^\dagger B^{\sss \top}\,X\,A
 +S\,R^\dagger B^{\sss \top}\,X\,B\,R^\dagger S^{\sss \top} \nonumber \\
\ns&\ns\ns&\ns  -A^{\sss \top} X\,B\,R_{\scriptscriptstyle X}^\dagger B^{\sss \top} X\, A
 +A^{\sss \top}  X\,B\,R_{\scriptscriptstyle X}^\dagger (B^{\sss \top} X\,B+R-R)\,R^\dagger S^{\sss \top}  \nonumber\\
 \ns&\ns\ns&\ns
 +S\,R^\dagger (B^{\sss \top} X\,B+R-R)\,R_{\scriptscriptstyle X}^\dagger B^{\sss \top} X\,A \nonumber
  \\
\ns&\ns\ns&\ns
 -S\,R^\dagger (B^{\sss \top} X\,B+R-R)\,R_{\scriptscriptstyle X}^\dagger (B^{\sss \top} X\,B+R-R)\,R^\dagger S^{\sss \top} +Q-S\,R^\dagger S^{\sss \top}  \nonumber \\
 %
 %
   \ns&\ns =\ns&\ns A^{\sss \top} \,X\,A
 -A^{\sss \top} \,X\,B\,R^\dagger S^{\sss \top} -S\,R^\dagger B^{\sss \top}\,X\,A
 +S\,R^\dagger B^{\sss \top}\,X\,B\,R^\dagger S^{\sss \top} \nonumber \\
\ns&\ns\ns&\ns  -A^{\sss \top} X\,B\,R_{\scriptscriptstyle X}^\dagger B^{\sss \top} X\, A
 +A^{\sss \top}  X\,B\,R_{\scriptscriptstyle X}^\dagger R_{\scriptscriptstyle X} \,R^\dagger S^{\sss \top}-A^{\sss \top}  X\,B\,R_{\scriptscriptstyle X}^\dagger  S^{\sss \top}  \nonumber\\
 \ns&\ns\ns&\ns
 +S\,R^\dagger R_{\scriptscriptstyle X}\,R_{\scriptscriptstyle X}^\dagger B^{\sss \top} X\,A
 -S\,\,R_{\scriptscriptstyle X}^\dagger B^{\sss \top} X\,A-S\,R^\dagger R_{\scriptscriptstyle X}\,R^\dagger S^{\sss \top} \nonumber
  \\
\ns&\ns\ns&\ns
  +S\,R^\dagger R_{\scriptscriptstyle X}\,R_{\scriptscriptstyle X}^\dagger S^{\sss \top}
   +S\,R_{\scriptscriptstyle X}^\dagger R_{\scriptscriptstyle X}\,R^\dagger S^{\sss \top}
 -S\,R_{\scriptscriptstyle X}^\dagger  S^{\sss \top}+Q-S\,R^\dagger S^{\sss \top}. \label{catena}
  \eea
From $\ker R_{\sss X} \subseteq \ker S_{\sss X}$, it follows that there exists $K$ such that $S_{\sss X}=K\,R_{\sss X}$, which gives
\bea
\label{newfact}
S_{\sss X}\,R_{\sss X}^\dagger \,R_{\sss X}=K\,R_{\sss X}\,R_{\sss X}^\dagger \,R_{\sss X}=K\,R_{\sss X}=S_{\sss X}.
\eea
Using this identity and its transpose, we can develop the terms in the right hand-side of the last equality sign of (\ref{catena}) as
\beann
A^{\sss \top} X \,B\,R_{\scriptscriptstyle X}^\dagger R_{\scriptscriptstyle X}\,R^\dagger S^{\sss \top}+S\,R_{\scriptscriptstyle X}^\dagger R_{\scriptscriptstyle X}\,R^\dagger\,S^{\sss \top}=S_{\scriptscriptstyle X}\,R_{\scriptscriptstyle X}^\dagger R_{\scriptscriptstyle X}\,R^\dagger S^{\sss \top}=S_{\scriptscriptstyle X}\,R^\dagger S^{\sss \top},
\eeann
\beann
S\,R^\dagger R_{\scriptscriptstyle X}\,R_{\scriptscriptstyle X}^\dagger B^{\sss \top} X A+S\,R^\dagger R_{\scriptscriptstyle X}\,R_{\scriptscriptstyle X}^\dagger S^{\sss \top}=S\,R^\dagger R_{\scriptscriptstyle X}\,R_{\scriptscriptstyle X}^\dagger S_{\scriptscriptstyle X}^{\sss \top}=S\,R^\dagger S_{\scriptscriptstyle X}^{\sss \top}.
\eeann
and
\beann
S\,R^\dagger B^{\sss \top}\,X\,B\,R^\dagger S^{\sss \top}-S\,R^\dagger R_{\scriptscriptstyle X}\,R^\dagger S^{\sss \top}=-S\,R^\dagger \,R\,R^\dagger\,S^{\sss \top}=-S\,R^\dagger S^{\sss \top}.
\eeann
Using these new simplified expressions back into (\ref{catena}) gives
\beann
X \ns&\ns =\ns&\ns 
  -A^{\sss \top} X \,B\,R^\dagger -S\,R^\dagger S^{\sss \top}-S\,R^\dagger B^{\sss \top}\,X\,A-S\,R^\dagger S^{\sss \top}+S_{\scriptscriptstyle X}\,R^\dagger S^{\sss \top}-S\,R^\dagger S_{\scriptscriptstyle X}^{\sss \top} \\
  \ns&\ns =\ns&\ns 
A^{\sss \top} X\,A-A^{\sss \top} X\,B\,R_{\scriptscriptstyle X}^\dagger B^{\sss \top} X \,A-S\,R_{\scriptscriptstyle X}^\dagger B^{\sss \top} X\,A-A^{\sss \top} X \,B\,R_{\scriptscriptstyle X}^\dagger S^{\sss \top}-S\,R_{\scriptscriptstyle X}^\dagger S^{\sss \top}+Q \\
 \ns&\ns\ns&\ns -(A^{\sss \top} X \,B+S)\,R^\dagger S^{\sss \top}-S\,R^\dagger (B^{\sss \top}\,X\,A+ S^{\sss \top})+S_{\scriptscriptstyle X}\,R^\dagger S^{\sss \top}-S\,R^\dagger S_{\scriptscriptstyle X}^{\sss \top} \\
\ns&\ns =\ns&\ns A^{\sss \top} X A  -  (A^{\sss \top} X B  +  S)(R  +  B^{\sss \top} X B)^\dagger(B^{\sss \top} X A  +  S^{\sss \top})  +  Q,
\eeann
which is indeed (\ref{gdare}).
We conclude the proof of {\bf (i)} showing that (\ref{kercond}) is equivalent to (\ref{kercond0}). We write (\ref{kercond0}) as
\beann
\ker R_{\scriptscriptstyle X}\ns&\ns \subseteq \ns&\ns \ker (A_{\scriptscriptstyle 0}^{\sss \top}\,X\,B) \\
\ns&\ns = \ns&\ns \ker (A^{\sss \top}\,X\,B-S\,R^\dagger \,B^{\sss \top}\,X\,B) \\
\ns&\ns = \ns&\ns \ker [A^{\sss \top}\,X\,B-S\,R^\dagger \,(R+B^{\sss \top}\,X\,B-R)] \\
\ns&\ns = \ns&\ns \ker (A^{\sss \top}\,X\,B+S-S\,R^\dagger\,R_{\scriptscriptstyle X}),
\eeann
since $S\,R^\dagger\,R=S$ in view of the second point in Lemma \ref{lem1}. Suppose (\ref{kercond}) holds. 
Let $\omega\in \ker R_{\scriptscriptstyle X}$. Then $S_{\scriptscriptstyle X}\,\omega=(S+A^{\sss \top}\,X\,B)\,\omega=0$. Thus, we have also $(A^{\sss \top}\,X\,B+S-S\,R^\dagger\,R_{\scriptscriptstyle X})\,\omega=0$ since $\omega\in \ker R_{\scriptscriptstyle X}$. Conversely, suppose that (\ref{kercond0}) holds true, and take $\omega\in \ker R_{\scriptscriptstyle X}$. Then, $(A^{\sss \top}\,X\,B+S-S\,R^\dagger\,R_{\scriptscriptstyle X})\,\omega=0$ implies $(S+A^{\sss \top}\,X\,B)\,\omega=0$.

Let us now consider {\bf (ii)}. We first show that $(R_{\scriptscriptstyle X}^{\dagger} \,R_{\scriptscriptstyle X} -I_{\sss m})\,R^\dagger=0$. {To prove this fact | which is trivial in the case of the  standard DARE($\Sigma$) | } we use the inclusion $\ker R_{\scriptscriptstyle X} \subseteq \ker R$, which holds true for any symmetric solution $X$ of CGDARE($\Sigma$), see
\cite[Lemma 4.1]{Ferrante-N-12}. 
In a suitable basis of the input space, $R_{\scriptscriptstyle X}$ can be written as
$R_{\scriptscriptstyle X}=\bsmat R_{\sss X,1} & 0 \\[1mm] 0 & 0 \esmat$, where $R_{\sss X,1}$ is invertible; let $\mu$ be the order of $R_{\sss X,1}$. In this basis, $R$ is written as $R=\bsmat R_{\sss 1} & 0 \\[1mm] 0 & 0 \esmat$, where $R_{\sss 1}$ may or may not be singular, and we obtain
\bea
\label{id}
(R_{\scriptscriptstyle X}^{\dagger} \,R_{\scriptscriptstyle X} -I_{\sss m})\,R^\dagger \ns&\ns = \ns&\ns
\left(\bmat{cc} R_{\sss X,1}^{-1} & 0 \\ 0 & 0 \emat \bmat{cc} R_{\sss X,1} & 0 \\ 0 & 0 \emat -\bmat{cc} I_{\mu} & 0 \\ 0 & I_{\sss m -\mu} \emat \right)\bmat{cc} R_{1}^{\dagger} & 0 \\ 0 & 0 \emat=0.
\eea
Thus,
\beann
A_{\scriptscriptstyle 0X} =A_{\scriptscriptstyle 0}-B\, (R  +  B^{\sss \top} X B)^\dagger B^{\sss \top} X A_{\scriptscriptstyle 0} \ns&\ns =\ns&\ns 
(A-B\,R^\dagger S^{\sss \top})-B\, (R  +  B^{\sss \top} X B)^\dagger B^{\sss \top} X (A-B\,R^\dagger S^{\sss \top}) \\
\ns&\ns =\ns&\ns A-B\,R^\dagger S^{\sss \top}-B\,R_{\scriptscriptstyle X}^\dagger B^{\sss \top} X\,A+B\,R_{\scriptscriptstyle X}^\dagger (R+B^{\sss \top} X\,B-R)\,R^\dagger S^{\sss \top} \\
\ns&\ns =\ns&\ns A_{\scriptscriptstyle X}+B\,(R_{\scriptscriptstyle X}^{\dagger} \,R_{\scriptscriptstyle X} -I_{\sss m})\,R^\dagger\,S^{\sss \top} = A_{\scriptscriptstyle X}.
\eeann
To prove {\bf (iii)} it suffices to observe that $Q_{\scriptscriptstyle 0}$ is the generalized Schur complement of $R$ in $\Pi$. Since $\Pi$ is assumed to be positive semidefinite, then such is also $Q_{\scriptscriptstyle 0}$.
\endproof
\ \\[-0.8cm]

Another useful result is the following generalization of a classic property of DARE($\Sigma$).\\

\begin{lemma}
\label{lemmadd}
Let  $T\in \real^{\sss n \times n}$ be invertible. Let 
\bea
A_{\scriptscriptstyle T} \defi T^{\sss -1} A_{\scriptscriptstyle 0}\,T, \quad B_{\scriptscriptstyle T} \defi T^{\sss -1} B, \quad Q_{\scriptscriptstyle T} \defi T^{\sss -1} Q_{\scriptscriptstyle 0}\,T.
\eea
Let also  $\Pi_{\sss T} \defi \bsmat Q_{\sss T} && 0 \\[1mm] 0 && R \esmat$ and $\Sigma_{\sss T} \defi (A_{\sss T},B_{\sss T},\Pi_{\sss T})$. 
 Then, $X$ is a solution of CGDARE($\Sigma$) -- and therefore also of CGDARE($\Sigma_{\sss 0}$) -- if and only if $X_{\scriptscriptstyle T}=T^{\sss -1} X\,T$ is a solution of CGDARE($\Sigma_{\sss T}$)
\bea
X_{\scriptscriptstyle T}  \ns&\ns=\ns&\ns   A_{\scriptscriptstyle T}^{\sss \top} X_{\scriptscriptstyle T} \,A_{\scriptscriptstyle T}  -  A_{\scriptscriptstyle T}^{\sss \top} X_{\scriptscriptstyle T} \,B_{\scriptscriptstyle T} \,(R  +  B_{\scriptscriptstyle T}^{\sss \top} X_{\scriptscriptstyle T} B_{\scriptscriptstyle T})^\dagger B_{\scriptscriptstyle T}^{\sss \top} X_{\scriptscriptstyle T} \,A_{\scriptscriptstyle T} +  Q_{\scriptscriptstyle T} \label{gdaretranf}\\
 &   &  \qquad \ker (R+B_{\scriptscriptstyle T}^{\sss \top}\,X_{\scriptscriptstyle T}\,B_{\scriptscriptstyle T}) \subseteq \ker (A_{\scriptscriptstyle T}^{\sss \top}\,X_{\scriptscriptstyle T}\,B_{\scriptscriptstyle T})\label{gkercondtranf}
\eea
\end{lemma}
\proof
The equations obtained by multiplying (\ref{gdare0}) to the left by $T^{\sss -1}$ and to the right by $T$  coincides with (\ref{gdaretranf}) with $X_{\scriptscriptstyle T}\defi T^{\sss -1} X\,T$. Moreover, since $T$ is invertible,
$\ker (R+B^{\sss \top}\,X\,B)  \subseteq \ker (A_{\scriptscriptstyle 0}^{\sss \top}\,X\,B)$
is equivalent to $\ker (R+B^{\sss \top}\,X\,B)  \subseteq  \ker (T^{\sss -1}\,A_{\scriptscriptstyle 0}^{\sss \top}\,X\,B)$,
which 
is equivalent to (\ref{gkercondtranf}).
\endproof

\section{Main results}
\subsection{Reduction corresponding to a singular $A_{\scriptscriptstyle 0}$}
\label{secA0}

In this section, we present the first fundamental result of this paper, that can be exploited as a basis for an iterative procedure -- to be used whenever $A_{\scriptscriptstyle 0}$ is singular -- to the end of decomposing the set of solutions of CGDARE($\Sigma$) into a trivial part and a part given by the set of solutions of a reduced order CGDARE.\\[-3mm]

\begin{theorem}
\label{the1}
Let $\nu \defi \dim(\ker A_{\scriptscriptstyle 0})$.
Let
$U=[\begin{array}{cc} U_{\sss 1} & U_{\sss 2} \end{array} ]$
be  an orthonormal change of coordinates in $\real^n$, where $\ima U_{\sss 2} = \ker A_{\scriptscriptstyle 0}$. Let $A_{\scriptscriptstyle U} \defi U^{\sss \top} A_{\scriptscriptstyle 0}\,U= [\begin{array}{cc} \tilde{A} & 0_{\sss n \times \nu} \end{array} ]$
where $\tilde{A}=\bsmat A_{\sss 1} \\[1mm] A_{\sss 21} \esmat$  with $A_{\sss 1} \in \real^{(n-\nu) \times (n-\nu)}$ and $A_{\sss 21} \in \real^{\nu \times (n-\nu)}$. Let also $B_{\scriptscriptstyle U}=U^{\sss \top} B$ and $Q_{\scriptscriptstyle U}=U^{\sss \top} Q_{\sss 0}\,U$ be partitioned conformably, i.e.,
 $B_{\scriptscriptstyle U}=\bsmat B_{\sss 1} \\[1mm] B_{\sss 2} \esmat$ and $Q_{\scriptscriptstyle U}=\bsmat Q_{\sss 11} & Q_{\sss 12} \\[1mm] Q_{\sss 12}^{\sss \top} & Q_{\sss 22} \esmat$, with $B_{\sss 1} \in \real^{(n-\nu) \times m}$, $B_{\sss 2} \in \real^{\nu \times m}$, $Q_{\sss 11}  \in \real^{(n-\nu) \times (n-\nu)}$ and $Q_{\sss 22}  \in \real^{\nu \times \nu}$.  Finally, let $Q_{\sss 1}  \defi   \tilde{A}^{\sss \top} \,Q_{\scriptscriptstyle U}\, \tilde{A}$, $S_{\sss 1}  \defi   \tilde{A}^{\sss \top}\,Q_{\scriptscriptstyle U}\,B_{\scriptscriptstyle U}$ and $R_{\sss 1}  \defi  R+B_{\scriptscriptstyle U}^{\sss \top}\,Q_{\scriptscriptstyle U}\,B_{\scriptscriptstyle U}$.

\begin{enumerate}
\item Let $X$ be a solution of CGDARE($\Sigma$), and partition $X_{\scriptscriptstyle U} \defi U^{\sss \top} X\,U$ as
$X_{\scriptscriptstyle U}=\bsmat X_{\sss 11} & X_{\sss 12} \\[1mm] X_{\sss 12}^{\sss \top} & X_{\sss 22} \esmat$, 
with $X_{\sss 11}  \in \real^{\sss (n-\nu) \times (n-\nu)}$ and $X_{\sss 22}  \in \real^{\sss \nu \times \nu}$. Then,
\begin{description}
\item{\bf (i)} there hold
\beann
X_{\sss 12}=Q_{\sss 12} \quad \text{and} \quad X_{\sss 22}=Q_{\sss 22}
\eeann
\item{\bf (ii)} The new Popov matrix
$\Pi_{\sss 1} \defi \bsmat Q_{\sss 1} & S_{\sss 1} \\[1mm] S_{\sss 1}^{\sss \top} & R_{\sss 1} \esmat$ is positive semidefinite. 
\item{\bf (iii)}
 Let $\Sigma_{\sss 1} \defi (A_{\sss 1},B_{\sss 1},\Pi_{\sss 1})$. Then, $\Delta_{\sss 1}\defi X_{\sss 11}-Q_{\sss 11}$ satisfies CGDARE($\Sigma_{\sss 1}$)
\bea
\Delta_{\sss 1}  \ns&\ns =\ns&\ns   A_{\sss 1}^{\sss \top} \Delta_{\sss 1}  A_{\sss 1}  -   (A_{\sss 1}^{\sss \top} \Delta_{\sss 1} B_{\sss 1}+S_{\sss 1}) (R_{\sss 1}  +  B_{\sss 1}^{\sss \top} \Delta_{\sss 1} B_{\sss 1})^\dagger (B_{\sss 1}^{\sss \top} \Delta_{\sss 1}  A_{\sss 1}+S_{\sss 1}^{\sss \top}) +  Q_{\sss 1} \label{redgdare} \\
 &   &  \qquad \ker (R_{\sss 1}+B_{\sss 1}^{\sss \top}\,\Delta_{\sss 1}\,B_{\sss 1}) \subseteq \ker (S_{\sss 1}+A_{\sss 1}^{\sss \top}\,\Delta_{\sss 1}\,B_{\sss 1}). \label{redkercond}
 \eea
\end{description}
\item Conversely, if $\Delta_{\sss 1}$ is a solution of (\ref{redgdare}-\ref{redkercond}), then
\bea
\label{newdefofX}
X=U\,\bmat{cc} \Delta_{\sss 1}+Q_{\sss 11} & Q_{\sss 12} \\ Q_{\sss 12}^\top & Q_{\sss 22} \emat\,U^\top
\eea
is a solution of CGDARE($\Sigma$). 
\end{enumerate}
\end{theorem}
\proof We begin proving the first point. In view of Lemma \ref{lemmadd},
$X$ is a solution of CGDARE($\Sigma$) if and only if $X_{\scriptscriptstyle U}=U^{\sss \top} X\,U$ is a solution of CGDARE($\Sigma_{\sss U})$
\bea
X_{\scriptscriptstyle U}  \ns&\ns=\ns&\ns   A_{\scriptscriptstyle U}^{\sss \top} X_{\scriptscriptstyle U} \,A_{\scriptscriptstyle U}  
-  A_{\scriptscriptstyle U}^{\sss \top} X_{\scriptscriptstyle U} \,B_{\scriptscriptstyle U} \,(R  +  B_{\scriptscriptstyle U}^{\sss \top} 
X_{\scriptscriptstyle U} B_{\scriptscriptstyle U})^\dagger B_{\scriptscriptstyle U}^{\sss \top} X_{\scriptscriptstyle U} \,A_{\scriptscriptstyle U} +  Q_{\scriptscriptstyle U} \label{gdaretranf1u}\\
 &   &  \qquad \ker (R+B_{\scriptscriptstyle U}^{\sss \top}\,X_{\scriptscriptstyle U}\,B_{\scriptscriptstyle U}) \subseteq \ker (A_{\scriptscriptstyle U}^{\sss \top}\,X_{\scriptscriptstyle U}\,B_{\scriptscriptstyle U}),\label{gkercondtranf1u}
\eea
where $\Pi_{\scriptscriptstyle U}=\bsmat Q_{\scriptscriptstyle U} && 0 \\[1mm] 0 && R \esmat$ and $\Sigma_{\scriptscriptstyle U}=(A_{\scriptscriptstyle U},B_{\scriptscriptstyle U},\Pi_{\scriptscriptstyle U})$.
Multiplying (\ref{gdaretranf1u}) 
 to the left by $[\begin{array}{cc} 0 & I_{\sss \nu} \end{array} ]$ yields
\beann
[\begin{array}{cc}  \! 0   \! & I_{\sss \nu}   \!  \end{array} ]  \! \bmat{cc}   \! X_{\sss 11}  \!  &  \!  X_{\sss 12}   \! \\   \! X_{\sss 12}^{\sss \top}   \! &   \! X_{\sss 22}  \!  \emat 
\ns&\ns   \! =  \! \ns&\ns
[\begin{array}{cc}   \!  0   \! & I_{\sss \nu}   \!  \end{array} ]  \! \bmat{cc}  \!  A_{1}^{\sss \top}   \! &   \! A_{\sss 21}^{\sss \top}   \! \\   \! 0   \! &  \!  0  \!   \emat X_{\scriptscriptstyle U} \, A_{\scriptscriptstyle U} \\
\ns&\ns\ns&\ns \hspace{-0.8cm}-[\begin{array}{cc}  \!    0   \! & I_{\sss \nu}    \! \end{array} ]  \! \bmat{cc}  \!  A_{1}^{\sss \top}   \! &  \!  A_{\sss 21}^{\sss \top}   \! \\   \! 0   \! &  \!  0  \!   \emat X_{\scriptscriptstyle U} \,B_{\scriptscriptstyle U} \,(R  +  B_{\scriptscriptstyle U}^{\sss \top} \,X_{\scriptscriptstyle U} \,B_{\scriptscriptstyle U})^\dagger\, B_{\scriptscriptstyle U}^{\sss \top}\, X_{\scriptscriptstyle U} \,A_{\scriptscriptstyle U}  \! +  \! [\begin{array}{cc}  \!   0  \!  & I_{\sss \nu}    \! \end{array} ]   \! \bmat{cc}  \!  Q_{\sss 11}   \! &   \! Q_{\sss 12}  \!  \\  \!  Q_{\sss 12}^{\sss \top}   \! &  \!  Q_{\sss 22}   \! \emat\!,
\eeann 
which gives 
$[\begin{array}{cc} X_{\sss 12}^{\sss \top} & X_{\sss 22}  \end{array} ]=[\begin{array}{cc} Q_{\sss 12}^{\sss \top} & Q_{\sss 22}  \end{array} ]$. This proves the first statement. 
To prove {\bf (ii)} we observe that
\bea
\label{eqPi1}
\Pi_{\sss 1} = \bmat{cc} Q_{\sss 1} & S_{\sss 1} \\ S_{\sss 1}^{\sss \top} & R_{\sss 1} \emat=\bmat{c} \tilde{A}^{\sss \top} \\ B^{\sss \top} \emat Q_{\scriptscriptstyle 0} [\begin{array}{cc} \tilde{A} & B \end{array} ]+\bmat{cc} 0 & 0 \\ 0 & R \emat \ge 0,
\eea
since, as shown in Lemma \ref{lem21}, $Q_{\scriptscriptstyle 0} \ge 0$. 
We now prove {\bf (iii)}.
Substitution of $X_{\scriptscriptstyle U}=Q_{\scriptscriptstyle U}+\bsmat \Delta_{\sss 1} & 0 \\[1mm] 0 & 0 \esmat$ obtained in the proof of {\bf (i)} into (\ref{gdaretranf1u}) gives
\beann
\bmat{cc}   \!  \Delta_{\sss 1}  \! & \! 0   \!   \\    \!  0  \! & \! 0   \!  \emat \ns\! & \!\ns = \ns\! & \!\ns 
\bmat{cc}  \!    Q_{\sss 1}  \! & \! 0  \!    \\   \!   0  \! & \! 0   \!  \emat+ \bmat{cc}   \!   A_{\sss 1}^{\sss \top}\,\Delta_{\sss 1}\,A_{\sss 1}  \! & \! 0    \!  \\     \! 0  \! & \! 0   \!  \emat-\bmat{c}   \! S_{\sss 1}+A_{\sss 1}^{\sss \top} \Delta_{\sss 1}\,B_{\sss 1}  \!  \\   \! 0   \! \emat (R_{\sss 1}+B_{\sss 1}^{\sss \top}\,\Delta_{\sss 1}\,B_{\sss 1})^{\dagger} [ \begin{array}{cc} S_{\sss 1}^{\sss \top}+B_{\sss 1}^{\sss \top} \,\Delta_{\sss 1}\,A_{\sss 1} \! & \! 0 \end{array} ],
\eeann
which is equivalent to (\ref{redgdare}). We now prove that $\Delta_{\sss 1}$ satisfies $\ker (R_{\sss 1}+B_{\sss 1}^{\sss \top}\,\Delta_{\sss 1}\,B_{\sss 1}) \subseteq \ker (S_{\sss 1}+A_{\sss 1}^{\sss \top}\,\Delta_{\sss 1}\,B_{\sss 1})$.
Substitution of $X_{\scriptscriptstyle U}=Q_{\scriptscriptstyle U}+\bsmat \Delta_{\sss 1} & 0 \\[1mm] 0 & 0 \esmat$ into (\ref{gkercondtranf1u}) gives
\beann
\ker (R_{\sss 1}+B_{\sss 1}^{\sss \top}\,\Delta_{\sss 1}\,B_{\sss 1}) \ns&\ns \subseteq  \ns&\ns \ker \left( \bmat{c} \tilde{A}^{\sss \top} \\ 0 \emat \,Q_{\scriptscriptstyle U}\,B_{\scriptscriptstyle U}+\bmat{cc} A_{\sss 1}^{\sss \top}\,\Delta_{\sss 1}\,B_{\sss 1} \\ 0 \emat  \right)= \ker \bmat{c} S_{\sss 1}+A_{\sss 1}^{\sss \top}\,\Delta_{\sss 1}\,B_{\sss 1}\\ 0 \emat,
\eeann
which is equivalent to (\ref{redkercond}). 
We now prove the converse. Let $X$ be as in (\ref{newdefofX}). Substituting $X_{\sss U}=U^\top\,X\,U=
\bsmat \Delta_{\sss 1}+Q_{\sss 11} && Q_{\sss 12} \\[1mm] Q_{\sss 12}^\top && Q_{\sss 22} \esmat$ 
 into CGDARE($\Sigma_{\sss U}$) gives
\beann
\bmat{cc} \!\! \Delta_{\sss 1}+Q_{\sss 11} \!\! & \!\! Q_{\sss 12} \!\! \\  \!\! Q_{\sss 12}^\top \!\! & \!\! Q_{\sss 22} \!\! \emat \ns&\ns  = \ns&\ns 
\bmat{cc} \!\! A_{\sss 1}^\top \!\! & \!\! A_{\sss 21}^\top \!\! \\  \!\! 0 \!\! & \!\! 0 \!\! \emat \bmat{cc} \!\! \Delta_{\sss 1}+Q_{\sss 11} \!\! & \!\! Q_{\sss 12} \!\! \\  \!\!Q_{\sss 12}^\top \!\! & \!\! Q_{\sss 22} \!\! \emat \bmat{cc} \!\! A_{\sss 1} \!\! & \!\! 0 \\ A_{\sss 21} \!\! & \!\! 0 \!\! \emat \\
\hspace{-1cm} \ns&\ns   \ns&\ns  +\bmat{cc} \!\! A_{\sss 1}^\top \!\! & \!\! A_{\sss 21}^\top  \!\! \\  \!\! 0 \!\! & \!\! 0 \!\! \emat \!\! \bmat{cc} \!\! \Delta_{\sss 1}+Q_{\sss 11} \!\! & \!\! Q_{\sss 12}  \!\! \\  \!\! Q_{\sss 12}^\top \!\! & \!\! Q_{\sss 22} \!\! \emat \!\! \bmat{c}  \!\! B_{\sss 1} \!\!  \\  \!\! B_{\sss 2} \!\! \emat \!\! \left(R+[\begin{array}{cc} \!\! B_{\sss 1}^\top \!\! & \!\! B_{\sss 2}^\top  \!\! \end{array} ] \!\!\bmat{cc} \!\! \Delta_{\sss 1}+Q_{\sss 11} \!\! & \!\! Q_{\sss 12}  \!\! \\  \!\! Q_{\sss 12}^\top \!\! & \!\! Q_{\sss 22} \!\! \emat \!\!\bmat{c}  \!\! B_{\sss 1} \!\! \\  \!\!B_{\sss 2} \!\! \emat\right)^\dagger \\
\hspace{-1cm}  \ns&\ns   \ns&\ns \times 
\bmat{c}  \!\!B_{\sss 1} \!\! \\  \!\!B_{\sss 2} \!\! \emat\bmat{cc} \!\! \Delta_{\sss 1}+Q_{\sss 11} \!\! & \!\! Q_{\sss 12} \!\! \\  \!\!Q_{\sss 12}^\top \!\! & \!\! Q_{\sss 22} \!\! \emat \bmat{cc} \!\! A_{\sss 1} \!\! & \!\! 0 \!\! \\ \!\! A_{\sss 21} \!\! & \!\! 0 \!\! \emat+\bmat{cc} \!\!Q_{\sss 11} \!\! & \!\! Q_{\sss 12} \!\! \\ \!\! Q_{\sss 12}^\top \!\! & \!\! Q_{\sss 22} \!\! \emat
\eeann
Developing the products and recalling that we have defined $Q_{\sss 1}  =  \tilde{A}^{\sss \top} \,Q_{\scriptscriptstyle U}\, \tilde{A}$, $S_{\sss 1}  =  \tilde{A}^{\sss \top}\,Q_{\scriptscriptstyle U}\,B_{\scriptscriptstyle U}$ and $R_{\sss 1}  =  R+B_{\scriptscriptstyle U}^{\sss \top}\,Q_{\scriptscriptstyle U}\,B_{\scriptscriptstyle U}$ gives
\beann
\bmat{cc} \!\! \Delta_{\sss 1}+Q_{\sss 11} \!\! & \!\! Q_{\sss 12} \!\! \\  \!\! Q_{\sss 12}^\top \!\! & \!\! Q_{\sss 22} \!\! \emat \ns&\ns  = \ns&\ns  \bmat{cc} \!\! A_{\sss 1}^\top\,\Delta_{\sss 1}\,A_{\sss 1}+Q_{\sss 1}  \!\! & \!\! 0 \!\! \\  \!\! 0 \!\! & \!\! 0 \!\! \emat 
 - \bmat{cc} \!\! A_{\sss 1}^\top\,\Delta_{\sss 1}\,B_{\sss 1}+S_{\sss 1}   \!\! \\  \!\! 0 \!\! \emat (R_{\sss 1}+B_{\sss 1}^\top\,\Delta_{\sss 1}\,B_{\sss 1})^\dagger [\begin{array}{cc} B_{\sss 1}^\top\,\Delta_{\sss 1}\,A_{\sss 1}+S_{\sss 1}^\top & 0 \end{array} ] \\
\ns&\ns\ns&\ns +\bmat{cc} \!\! Q_{\sss 11} \!\! & \!\! Q_{\sss 12} \!\! \\  \!\! Q_{\sss 12}^\top \!\! & \!\! Q_{\sss 22} \!\! \emat,
 \eeann
 which is satisfied since $\Delta_{\sss 1}$ is a solution of (\ref{redgdare}-\ref{redkercond}).
\endproof

The following property, which considers the structure of the closed-loop matrix in the basis described by $U$, is stated separately from properties {\bf (i-iii)} in Theorem \ref{the1} to emphasize the differences between this first reduction and the second reduction that will be  presented in the next section. In fact, while in the standard case of DARE($\Sigma$)  this property of the closed-loop matrix applies to both the first and the second reduction procedure, in the general case of CGDARE($\Sigma$) the structure of the closed-loop matrix described in the following property is maintained only for the first reduction procedure. \\[-3mm]

\begin{proposition}
Given a solution $X$ of CGDARE($\Sigma$) and the associated solution $\Delta_{\sss 1}$ of (\ref{redgdare}-\ref{redkercond}), let $A_{\scriptscriptstyle X}$ and $A_{\sss \Delta_{\sss 1}}$ be the associated closed-loop matrices. Then,
\beann
U^{\sss \top} A_{\scriptscriptstyle X}\,U=\bmat{ll} A_{\sss \Delta_{\sss 1}} & 0 \\ \star & 0_{\sss \nu \times \nu} \emat.
\eeann
\end{proposition}
\proof
We first observe that the last $\nu$ columns of $U^{\sss \top} A_{\scriptscriptstyle X}\,U$ are also zero, i.e.,
\beann
U^{\sss \top} A_{\scriptscriptstyle X}\,U \ns&\ns = \ns&\ns U^{\sss \top} (A_{\scriptscriptstyle 0}-B\,R_{\scriptscriptstyle X}^\dagger B^{\sss \top} \,X\,A_{\scriptscriptstyle 0})U \\
\ns&\ns = \ns&\ns A_{\scriptscriptstyle U}-B_{\scriptscriptstyle U}\,(R+B_{\scriptscriptstyle U}^{\sss \top}\,X_{\scriptscriptstyle U}\,B_{\scriptscriptstyle U})^{\dagger}\,B_{\scriptscriptstyle U}^{\sss \top}\,X_{\scriptscriptstyle U}\,A_{\scriptscriptstyle U}=[ \begin{array}{cc}  \star & 0 \end{array} ],
\eeann
in view of the fact that the last $\nu$ columns of $A_{\scriptscriptstyle U}$ are zero. Moreover,
\beann
U^{\sss \top} A_{\scriptscriptstyle X}\,U \ns&\ns = \ns&\ns \bmat{cc} A_{\sss 1} & 0 \\ A_{\sss 21} & 0 \emat - \bmat{c} B_{\sss 1}\\ B_{\sss 2}\emat \left[ R+ [ \begin{array}{cc} B_{\sss 1}^{\sss \top} & B_{\sss 2}^{\sss \top} \end{array}]\left(Q_{\scriptscriptstyle U}+\bmat{cc} \Delta_{\sss 1} & 0 \\ 0 & 0 \emat \right)\bmat{c} B_{\sss 1}\\ B_{\sss 2} \emat \right]^\dagger B_{\scriptscriptstyle U}\,X_{\scriptscriptstyle U}\,A_{\scriptscriptstyle U} \\
 \ns&\ns = \ns&\ns \bmat{cc} A_{\sss 1} & 0 \\ A_{\sss 21} & 0 \emat - \bmat{c} B_{\sss 1}\\ B_{\sss 2}\emat (R_{\sss 1}+ B_{\sss 1}^{\sss \top}  \Delta_{\sss 1}\,B_{\sss 1})^\dagger B_{\scriptscriptstyle U}^{\sss \top}\,X_{\scriptscriptstyle U}\,A_{\scriptscriptstyle U}
\eeann
and
\beann
A_{\Delta_{\sss 1}} \ns&\ns = \ns&\ns A_{\sss 1}-B_{\sss 1}\,(R_{\sss 1}+B_{\sss 1}^{\sss \top} \Delta_{\sss 1}\,B_{\sss 1})^\dagger (B_{\sss 1}^{\sss \top} \Delta_{\sss 1}\,A_{\sss 1}+S_{\sss 1}^{\sss \top})-B_{\sss 1}\,R_{\sss 1}^\dagger S_{\sss 1}^{\sss \top} +B_{\sss 1}\,R_{\sss 1}^\dagger S_{\sss 1}^{\sss \top} \\
\ns&\ns = \ns&\ns A_{\sss 1}-B_{\sss 1}\,(R_{\sss 1}+B_{\sss 1}^{\sss \top} \Delta_{\sss 1}\,B_{\sss 1})^\dagger B_{\sss 1}^{\sss \top} \Delta_{\sss 1}\,A_{\sss 1}-B_{\sss 1}\,(R_{\sss 1}+B_{\sss 1}^{\sss \top} \Delta_{\sss 1}\,B_{\sss 1})^\dagger R_{\sss 1}\,R_{\sss 1}^\dagger\,S_{\sss 1}^{\sss \top} \\
\ns&\ns\ns&\ns
-B_{\sss 1}\,R_{\sss 1}^\dagger S_{\sss 1}^{\sss \top} +B_{\sss 1}\,(R_{\sss 1}+B_{\sss 1}^{\sss \top} \Delta_{\sss 1}\,B_{\sss 1})^\dagger(R_{\sss 1}+B_{\sss 1}^{\sss \top} \Delta_{\sss 1}\,B_{\sss 1}) R_{\sss 1}^\dagger S_{\sss 1}^{\sss \top},
\eeann
where the last equality follows from the identity $(R_{\sss 1}+B_{\sss 1}^{\sss \top} \Delta_{\sss 1}\,B_{\sss 1})^\dagger(R_{\sss 1}+B_{\sss 1}^{\sss \top} \Delta_{\sss 1}\,B_{\sss 1}) R_{\sss 1}^\dagger=R_{\sss 1}^\dagger$, which can be proved exactly in the same way as (\ref{id}). \footnote{Indeed, in CGDARE($\Sigma_{\sss 1}$) the matrices $R_{\sss 1}$ and $R_{\sss 1}+B_{\sss 1}^{\sss \top}\,\Delta_{\sss 1}\,B_{\sss 1}$ play the same role of $R$ and $R+B^{\sss \top}\,X\,B$ in CGDARE($\Sigma$), so that $\ker (R_{\sss 1}+B_{\sss 1}^{\sss \top}\,\Delta_{\sss 1}\,B_{\sss 1})\subseteq \ker R_{\sss 1}$.}
Thus,
\beann
A_{\Delta_{\sss 1}} \ns&\ns = \ns&\ns A_{\sss 1}-B_{\sss 1}\,(R_{\sss 1}+B_{\sss 1}^{\sss \top} \Delta_{\sss 1}\,B_{\sss 1})^\dagger B_{\sss 1}^{\sss \top} \Delta_{\sss 1}\,A_{\sss 1}-B\,R_{\sss 1}^\dagger S_{\sss 1}^{\sss \top} 
-B_{\sss 1}\,(R_{\sss 1}+B_{\sss 1}^{\sss \top} \Delta_{\sss 1}\,B_{\sss 1})^\dagger (R_{\sss 1}-R_{\sss 1}-B_{\sss 1}^{\sss \top} \Delta_{\sss 1}\,B_{\sss 1}) R_{\sss 1}^\dagger S_{\sss 1}^{\sss \top} \\
\ns&\ns = \ns&\ns A_{\sss 1}-B\,R_{\sss 1}^\dagger S_{\sss 1}^{\sss \top} -B_{\sss 1}\,(R_{\sss 1}+B_{\sss 1}^{\sss \top} \Delta_{\sss 1}\,B_{\sss 1})^\dagger B_{\sss 1}^{\sss \top} \Delta_{\sss 1}\,A_{\sss 1}
+B_{\sss 1}\,(R_{\sss 1}+B_{\sss 1}^{\sss \top} \Delta_{\sss 1}\,B_{\sss 1})^\dagger B_{\sss 1}^{\sss \top} \Delta_{\sss 1}\,B_{\sss 1}\, R_{\sss 1}^\dagger S_{\sss 1}^{\sss \top} \\
\ns&\ns = \ns&\ns A_{\sss 1}-B\,R_{\sss 1}^\dagger S_{\sss 1}^{\sss \top} -B_{\sss 1}\,(R_{\sss 1}+B_{\sss 1}^{\sss \top} \Delta_{\sss 1}\,B_{\sss 1})^\dagger B_{\sss 1}^{\sss \top} \Delta_{\sss 1}\,(A_{\sss 1}-B_{\sss 1}\, R_{\sss 1}^\dagger S_{\sss 1}^{\sss \top}).
\eeann
Then, denoting by $\Gamma$ the upper-left block submatrix of order $n-\nu$ within $U^{\sss \top} A_{\scriptscriptstyle X}\,U$, we find
\bea
\Gamma-A_{\Delta_{\sss 1}} \ns&\ns = \ns&\ns   B_{\sss 1}\,(R_{\sss 1}+B_{\sss 1}^{\sss \top} \Delta_{\sss 1}\,B_{\sss 1})^\dagger(B_{\sss 1}^{\sss \top}\,\Delta_{\sss 1}\,A_{\sss 1}-B_{\scriptscriptstyle U}^{\sss \top}\,X_{\scriptscriptstyle U}\,\tilde{A}) \nonumber\\
\ns&\ns\ns&\ns+B_{\sss 1}\,R_{\sss 1}^\dagger\,S_{\sss 1}^{\sss \top}-B_{\sss 1}\,(R_{\sss 1}+B_{\sss 1}^{\sss \top} \Delta_{\sss 1}\,B_{\sss 1})^\dagger\,B_{\sss 1}^{\sss \top}\,\Delta_{\sss 1}\,B_{\sss 1}\,R_{\sss 1}^\dagger \,S_{\sss 1}^{\sss \top}. \label{aux}
\eea
A simple calculation shows also that
\beann
B_{\sss 1}^{\sss \top}\,\Delta_{\sss 1}\,A_{\sss 1}-B_{\scriptscriptstyle U}^{\sss \top}\,X_{\scriptscriptstyle U}\,\tilde{A}  \ns&\ns = \ns&\ns B_{\sss 1}^{\sss \top}\,\Delta_{\sss 1}\,A_{\sss 1}-[\begin{array}{cc} B_{\sss 1}^{\sss \top} & B_{\sss 2}^{\sss \top} \end{array} ] \bmat{cc} (Q_{\sss 11}+\Delta_{\sss 1})& Q_{\sss 12} \\ Q_{\sss 12}^{\sss \top} & Q_{\sss 22} \emat \bmat{c} A_{1} \\ A_{\sss 21} \emat \\
 \ns&\ns = \ns&\ns -[\begin{array}{cc} B_{\sss 1}^{\sss \top} & B_{\sss 2}^{\sss \top} \end{array} ]  \bmat{cc} Q_{\sss 11}& Q_{\sss 12} \\ Q_{\sss 12}^{\sss \top} & Q_{\sss 22} \emat \bmat{c} A_{1} \\ A_{\sss 21} \emat -B_{\scriptscriptstyle U}^{\sss \top}\,Q_{\scriptscriptstyle U}\,\tilde{A}=-S_{\sss 1}^{\sss \top}.
\eeann
We can use this identity in (\ref{aux}) and we obtain
\beann
\Gamma -A_{\Delta_{\sss 1}} \ns&\ns = \ns&\ns -B_{\sss 1}\,(R_{\sss 1}+B_{\sss 1}^{\sss \top} \Delta_{\sss 1}\,B_{\sss 1})^\dagger\,S_{\sss 1}^{\sss \top}+B_{\sss 1}\,R_{\sss 1}^\dagger S_{\sss 1}^{\sss \top}-B_{\sss 1}\,(R_{\sss 1}+B_{\sss 1}^{\sss \top} \Delta_{\sss 1}\,B_{\sss 1})^\dagger B_{\sss 1}^{\sss \top} \Delta_{\sss 1}\,B_{\sss 1}\,R_{\sss 1}^\dagger S_{\sss 1}^{\sss \top} \\
\ns&\ns = \ns&\ns -B_{\sss 1}\,(R_{\sss 1}+B_{\sss 1}^{\sss \top} \Delta_{\sss 1}\,B_{\sss 1})^\dagger\,S_{\sss 1}^{\sss \top}+B_{\sss 1}\,(R_{\sss 1}+B_{\sss 1}^{\sss \top} \Delta_{\sss 1}\,B_{\sss 1})^\dagger(R_{\sss 1}+B_{\sss 1}^{\sss \top} \Delta_{\sss 1}\,B_{\sss 1})\,R_{\sss 1}^\dagger S_{\sss 1}^{\sss \top}\\
\ns&\ns \ns&\ns -B_{\sss 1}\,(R_{\sss 1}+B_{\sss 1}^{\sss \top} \Delta_{\sss 1}\,B_{\sss 1})^\dagger B_{\sss 1}^{\sss \top} \Delta_{\sss 1}\,B_{\sss 1}\,R_{\sss 1}^\dagger S_{\sss 1}^{\sss \top}\\
\ns&\ns = \ns&\ns B_{\sss 1}\,(R_{\sss 1}+B_{\sss 1}^{\sss \top} \Delta_{\sss 1}\,B_{\sss 1})^\dagger \left[(R_{\sss 1}+B_{\sss 1}^{\sss \top} \Delta_{\sss 1}\,B_{\sss 1})\,R_{\sss 1}^\dagger \,S_{\sss 1}^{\sss \top}-S_{\sss 1}^{\sss \top}-B_{\sss 1}^{\sss \top} \Delta_{\sss 1}\,B_{\sss 1}\,R_{\sss 1}^\dagger S_{\sss 1}^{\sss \top} \right]\\
\ns&\ns = \ns&\ns B_{\sss 1}\,(R_{\sss 1}+B_{\sss 1}^{\sss \top} \Delta_{\sss 1}\,B_{\sss 1})^\dagger (R_{\sss 1}\,R_{\sss 1}^\dagger \,S_{\sss 1}^{\sss \top}-S_{\sss 1}^{\sss \top} )=0.
\eeann
\endproof

In view of {\bf (i)} of Theorem \ref{the1}, all solutions  of CGDARE($\Sigma$) coincide along the subspace
$\gU \defi \ker \left( \bsmat I_{n-\nu} & 0 \\[1mm] 0 & 0 \esmat \,U^{\sss \top}\right)$. This means that given any two solutions $X$ and $Y$ of CGDARE($\Sigma$), we have $X|_{\gU}=Y|_{\gU}=Q_{\scriptscriptstyle 0}|_{\gU}$.

The following result gives a property of the set of solutions of CGDARE($\Sigma$), and a procedure to solve CGDARE($\Sigma$) in terms of the reduced order DARE($\Sigma$).

\begin{corollary}
\label{cor1}
The set $\gX$ of solutions of CGDARE($\Sigma$) is parameterized as the set of matrices that can be expressed as
\[
X=U\,\bsmat \Delta_{\sss 1} && 0 \\[1mm] 0 && 0 \esmat \,U^{\sss \top}+Q_{\scriptscriptstyle 0}
\]
where $U=[\begin{array}{cc} U_{\sss 1} & U_{\sss 2} \end{array} ]$ is defined as in Theorem \ref{the1} and $\Delta_{\sss 1}$ is solution of (\ref{redgdare}-\ref{redkercond}).\\
\end{corollary}

{
After the reduction described in Theorem \ref{the1}, it may still happen that $A_{\sss 1} -B_{\sss 1}\,R_{\sss 1}^\dagger\,S_{\sss 1}$ is singular. 
However, since we have proved that CGDARE($\Sigma_{\sss 1}$) has exactly the same structure of CGDARE($\Sigma$), because $\Pi_{\sss 1} =\Pi_{\sss 1}^{\sss \top} \ge 0$, if $A_{\sss 1} -B_{\sss 1}\,R_{\sss 1}^\dagger\,S_{\sss 1}$ is singular we can iterate the procedure by rewriting (\ref{redgdare}-\ref{redkercond}) as
\bea
\label{redgdare1}
\Delta_{\sss 1}  \ns&\ns =\ns&\ns   A_{\sss 0,1}^{\sss \top} \Delta_{\sss 1}  A_{\sss 0,1}  -   A_{\sss 0,1}^{\sss \top} \Delta_{\sss 1} B_{\sss 1} (R_{\sss 1}  +  B_{\sss 1}^{\sss \top} \Delta_{\sss 1} B_{\sss 1})^\dagger B_{\sss 1}^{\sss \top} \Delta_{\sss 1}  A_{\sss 0,1} +  Q_{\sss 0,1} 
 \\
 &   &  \qquad \ker (R_{\sss 1}+B_{\sss 1}^{\sss \top}\,\Delta_{\sss 1}\,B_{\sss 1}) \subseteq \ker (A_{\sss 0,1}^{\sss \top}\,\Delta\,B_{\sss 1}), \label{redkercond1}
 \eea
where $A_{\sss 0,1} \defi A_{\sss 1} -B_{\sss 1}\,R_{\sss 1}^\dagger\,S_{\sss 1}^{\sss \top}$  and $Q_{\sss 0,1} \defi Q_{\sss 1} -S_{\sss 1}\,R_{\sss 1}^\dagger\,S_{\sss 1}^{\sss \top}$, and choosing a basis where $A_{\sss 0,1}=[\begin{array}{cc} \tilde{A}_{\sss 1} & 0 \end{array} ]$ and $\tilde{A}_{\sss 1}$ is of full column-rank. 
By following iteratively the procedure that led from CGDARE($\Sigma$) to CGDARE($\Sigma_{\sss 1}$), we eventually obtain a CGDARE($\Sigma_k$) of the form 
\bea
\Delta_{\sss k}  \ns&\ns =\ns&\ns   A_{\sss 0,k}^{\sss \top} \Delta_{\sss k}\,  A_{\sss 0,k}  -   A_{\sss 0,k}^{\sss \top} \Delta_{\sss k}\,B_{\sss k}\, (R_{\sss k}  +  B_{\sss k}^{\sss \top} \Delta_{\sss k}\, B_{\sss k})^\dagger\, B_{\sss k}^{\sss \top} \Delta_{\sss k}\,  A_{\sss 0,k} +  Q_{\sss 0,k} \label{redgdarek} \\
 &   &  \qquad \ker (R_{\sss k}+B_{\sss k}^{\sss \top}\,\Delta_{\sss k}\,B_{\sss k}) \subseteq \ker (A_{\sss 0,k}^{\sss \top}\,\Delta_{\sss k}\,B_{\sss k}), \label{redkercondk}
 \eea
 where now $A_{\sss 0,k}$ is non-singular. Notice also that  this reduction procedure can be carried out 
 only using the problem data $A,B,Q,R,S$, so that it holds for any solution $X$ of CGDARE($\Sigma$). In other words, this procedure (and the one that will follow in the next section) can be performed without the need to compute a particular solution of the Riccati equation.
  
 Once we have obtained the reduced-order CGDARE, if the corresponding matrix $R$ is singular, we can proceed with the second reduction procedure outlined in the next section.
}

\subsection{Reduction corresponding to a singular $R$}
\label{sec2b}
Consider CGDARE($\Sigma)$, either in the form given by (\ref{gdare}-\ref{kercond}) or (\ref{gdare0}-\ref{kercond0}). Suppose $R$ is singular. We assume that we have already performed the reduction described in the previous section. Hence, we may assume that $A_{\scriptscriptstyle 0}$ is now non-singular. 
To deal with this situation, we address separately two different cases: the first leads either to a reduced-order DARE or to a symmetric Stein equation depending on the rank of $R$, and the second leads to a reduced-order CGDARE.
We first consider the case in which $A_{\scriptscriptstyle 0}^{-1} B\,\ker R = \{0\}$, i.e., $B\,\ker R = \{0\}$. This case can in turn be divided into two sub-cases. The first is the one in which $R$ is not the zero matrix.
In this case, denoting by $r$ the rank of $R$, we can consider a change of coordinates in the input space  that brings $R$ in the form
\[
R=\bmat{cc} R_{\sss 1} & 0 \\ 0 & 0 \emat,
\]
where $R_{\sss 1}$ is non-singular, and $r$ is its order. With respect to this basis, since $\ker R=\ima \bsmat 0 \\[1mm] I_{m-r} \esmat$, matrix $B$ can be written as $B=[\begin{array}{cc} B_{\sss 1} & 0_{n \times (m-r)} \end{array} ]$, and (\ref{gdare0}-\ref{kercond0}) written in this basis

\beann
X \ns&\ns = \ns&\ns A_{\sss 0}^{\sss \top}\,X\,A_{\sss 0}-A_{\sss 0}^{\sss \top}\,[\begin{array}{cc} X\,B_{\sss 1} & 0 \end{array}]\left( 
\bmat{cc} R_1 & 0 \\ 0 & 0 \emat+\bmat{cc}B_{\sss 1}^{\sss \top} X\,B_{\sss 1}  & 0 \\ 0 & 0 \emat\right)^\dagger
[ \begin{array}{cc} B_{\sss 1}^{\sss \top} X & 0 \end{array} ] \,A_{\sss 0}+Q_{\sss 0}\\
&& \ker  \left( 
\bmat{cc} R_1 +B_{\sss 1}^{\sss \top} X\,B_{\sss 1}  & 0 \\ 0 & 0 \emat\right) \subseteq \ker (A_{\sss 0}^{\sss \top} X\,[ \begin{array}{cc} B_{\sss 1} & 0  \end{array} ]),
\eeann
reduces to
\beann
X \ns&\ns = \ns&\ns A_{\sss 0}^{\sss \top}\,X\,A_{\sss 0}-A_{\sss 0}^{\sss \top}\,X\,B_{\sss 1} (R_1 +B_{\sss 1}^{\sss \top} X\,B_{\sss 1})^\dagger
B_{\sss 1}^{\sss \top} X \, A_{\sss 0}+Q_{\sss 0} \\
&& \ima \bmat{c} 0 \\ I_{\sss m-r} \emat \subseteq \ker [ \begin{array}{cc} \star & 0_{\sss n \times (m-r)} \end{array} ]
\eeann
where now $R_1$ is invertible as required, so that $R_1 +B_{\sss 1}^{\sss \top} X\,B_{\sss 1}$ is positive definite. Hence, the latter is in fact a DARE
\beann
X \ns&\ns = \ns&\ns A_{\sss 0}^{\sss \top}\,X\,A_{\sss 0}-A_{\sss 0}^{\sss \top}\,X\,B_{\sss 1} (R_1 +B_{\sss 1}^{\sss \top} X\,B_{\sss 1})^{-1}
B_{\sss 1}^{\sss \top} X \, A_{\sss 0}+Q_{\sss 0}.
\eeann
 If $r=0$, i.e., if $R$ is the zero matrix, then $B\,\ker R=\{0\}$ implies that $B$ is also the zero matrix. In this case, CGDARE($\Sigma$) reduces to a symmetric Stein equation\footnote{For a discussion on the properties of symmetric Stein equations we refer to \cite[Section 5.3]{Lancaster-95} and \cite[Section 1.5]{Ionescu-OW-99}.}
\beann
X \ns&\ns = \ns&\ns A_{\sss 0}^{\sss \top}\,X\,A_{\sss 0}+Q_{\sss 0}.
\eeann

We now consider the case in which $A_{\scriptscriptstyle 0}^{-1} B\,\ker R \neq \{0\}$.

\begin{theorem}
\label{the2}
Let $\eta \defi \dim (A_{\scriptscriptstyle 0}^{-1} B\,\ker R)$.  Let
$V=[\begin{array}{cc} V_{\sss 1} & V_{\sss 2} \end{array} ]$
be an orthonormal change of coordinates in $\real^n$ where $\ima V_{\sss 2}=A_{\scriptscriptstyle 0}^{-1} B\,\ker R$. 
Let  $Q_{\scriptscriptstyle V} \defi V^{\sss \top}\,A_{\scriptscriptstyle 0}\,V$ and
$A_{\scriptscriptstyle V} \defi V^{\sss \top} A_{\scriptscriptstyle 0}\,V=\bsmat A_{\sss 1} & \star \\[1mm] \star & \star \esmat$, $B_{\scriptscriptstyle V}\defi V^{\sss \top}\,B=\bsmat B_{\sss 1} \\[1mm] \star \esmat$, $R_{\sss 1}\defi R+B^{\sss \top}\,Q_{\scriptscriptstyle 0}\,B$,
with $A_{\sss 1} \defi V_{\sss 1}^{\sss \top} A_{\scriptscriptstyle 0}\,V_{\sss 1} \in \real^{(n-\eta) \times (n-\eta)}$ and $B_{\sss 1}\defi V_{\sss 1}^{\sss \top} B \in \real^{(n-\eta) \times m}$. 
Let $Q_{\scriptscriptstyle V} \defi V^{\sss \top} Q_{\scriptscriptstyle 0}\,V= \bsmat Q_{\sss 11} & Q_{\sss 12} \\[1mm] Q_{\sss 12}^{\sss \top} & Q_{\sss 22} \esmat$, $A_{\scriptscriptstyle V}^{\sss \top}\,Q_{\scriptscriptstyle V}\,A_{\scriptscriptstyle V}=\bsmat Q_{\sss 1} & \star
\\[1mm] \star & \star \esmat$, $ A_{\scriptscriptstyle V}^{\sss \top} Q_{\scriptscriptstyle V}\,B_{\scriptscriptstyle V}=\bsmat S_{\sss 1} \\[1mm] \star \esmat$, 
where $Q_{\sss 11},Q_{\sss 1}\in \real^{(n-\eta) \times (n-\eta)}$ and $S_{\sss 1}\in \real^{(n-\eta) \times m}$. 
 Then, 
 \begin{enumerate}
 \item Let $X$ be a solution of CGDARE($\Sigma$), and partition $X_{\scriptscriptstyle V} \defi V^{\sss \top} X\,V$ as $X_{\scriptscriptstyle V} =\bsmat X_{\sss 11} & X_{\sss 12} \\[1mm] X_{\sss 12}^{\sss \top} & X_{\sss 22} \esmat$. Then,  
\begin{description}
\item{\bf (i)} there hold
\beann
X_{\sss 12}=Q_{\sss 12} \quad \text{and} \quad X_{\sss 22}=Q_{\sss 22}
\eeann
\item{\bf (ii)} The Popov matrix $\Pi_{\sss 1} \defi \bsmat Q_{\sss 1}& S_{\sss 1} \\[1mm] S_{\sss 1}^{\sss \top} & R_{\sss 1}\esmat$ is positive semidefinite. 
\item{\bf (iii)} Let $\Sigma_{\sss 1} \defi (A_{\sss 1},B_{\sss 1},\Pi_{\sss 1})$. Then, $\Delta_{\sss 1} \defi X_{\sss 11}-Q_{\sss 11}$ satisfies CGDARE($\Sigma_{\sss 1}$)
\bea
\Delta_{\sss 1}  \ns&\ns =\ns&\ns   A_{\sss 1}^{\sss \top} \Delta_{\sss 1}  A_{\sss 1}  -   (A_{\sss 1}^{\sss \top} \Delta_{\sss 1} B_{\sss 1}+S_{\sss 1}) (R_{\sss 1}  +  B_{\sss 1}^{\sss \top} \Delta_{\sss 1} B_{\sss 1})^\dagger (B_{\sss 1}^{\sss \top} \Delta_{\sss 1}  A_{\sss 1}+S_{\sss 1}^{\sss \top}) +  Q_{\sss 1} \label{redgdare2} \\
 &   &  \qquad \ker (R_{\sss 1}+B_{\sss 1}^{\sss \top}\,\Delta_{\sss 1}\,B_{\sss 1}) \subseteq \ker (S_{\sss 1}+A_{\sss 1}^{\sss \top}\,\Delta\,B_{\sss 1}). \label{redkercond2}
 \eea
\end{description}
\item Conversely, if $\Delta_{\sss 1}$ is a solution of (\ref{redgdare2}-\ref{redkercond2}), then
\[
X=V\,\bmat{cc} \Delta_{\sss 1}+Q_{\sss 11} & Q_{\sss 12} \\  Q_{\sss 12}^\top  & Q_{\sss 22} \emat\,V^\top
\] 
is a solution of CGDARE($\Sigma$). 
\end{enumerate}
\end{theorem}
\proof We prove the first point. As already observed in the beginning of Section \ref{secA0}, $X$ is a solution of (\ref{gdare}-\ref{kercond}) -- and therefore also of (\ref{gdare0}-\ref{kercond0}) -- if and only if $X_{\scriptscriptstyle V}=V^{\sss \top} X\,V$ is a solution of CGDARE($\Sigma_{\sss V})$
\bea
X_{\scriptscriptstyle V}  \ns&\ns=\ns&\ns   A_{\scriptscriptstyle V}^{\sss \top} X_{\scriptscriptstyle V} \,A_{\scriptscriptstyle V}  -  A_{\scriptscriptstyle V}^{\sss \top} X_{\scriptscriptstyle V} \,B_{\scriptscriptstyle V} \,(R  +  B_{\scriptscriptstyle V}^{\sss \top} X_{\scriptscriptstyle V} B_{\scriptscriptstyle V})^\dagger B_{\scriptscriptstyle V}^{\sss \top} X_{\scriptscriptstyle V} \,A_{\scriptscriptstyle V} +  Q_{\scriptscriptstyle V} \label{gdaretranf1}\\
 &   &  \qquad \ker (R+B_{\scriptscriptstyle V}^{\sss \top}\,X_{\scriptscriptstyle V}\,B_{\scriptscriptstyle V}) \subseteq \ker (A_{\scriptscriptstyle V}^{\sss \top}\,X_{\scriptscriptstyle V}\,B_{\scriptscriptstyle V}),\label{gkercondtranf1}
\eea
where $\Pi_{\scriptscriptstyle V}=\bsmat Q_{\scriptscriptstyle V} && 0 \\[1mm] 0 && R \esmat$ and $\Sigma_{\scriptscriptstyle V}=(A_{\scriptscriptstyle V},B_{\scriptscriptstyle V},\Pi_{\scriptscriptstyle V})$.
We can re-write (\ref{gdaretranf1}) as
\beann
X_{\sss V} \ns&\ns=\ns&\ns  
A_{\scriptscriptstyle V}^{\sss \top} X_{\scriptscriptstyle V} \,V^{\sss \top}\,[\,I_{\sss n}-
 B\,(R+B^{\sss \top}\,X\,B)^\dagger B^{\sss \top} X]\,A\,V+Q_{\sss V}.
\eeann
Post-multiplying the latter by $\bsmat 0 \\[1mm] I_{\eta}\esmat$ and considering a basis matrix $K_{\sss R}$ for $\ker R$, so that we can write $V_{\sss 2}=A^{-1} B\,K_{\sss R}$, gives
\beann
\bmat{c} X_{\sss 12} \\ X_{\sss 22}
\emat \ns&\ns=\ns&\ns A_{\scriptscriptstyle V}^{\sss \top} X_{\scriptscriptstyle V} \,V^{\sss \top}\,[\,I_{\sss n}-
 B\,(R+B^{\sss \top}\,X\,B)^\dagger B^{\sss \top} X]\,A\,V_{\sss 2}+\bmat{c} Q_{\sss 12} \\ Q_{\sss 22}
\emat \\
\ns&\ns=\ns&\ns V^{\sss \top}\,A_{\scriptscriptstyle 0}^{\sss \top}\,X\,B\,[\,I_{\sss m}-
 R_{\sss X}^\dagger (B^{\sss \top} X\,B+R-R)]\,K_{\sss R}+\bmat{c} Q_{\sss 12} \\ Q_{\sss 22}
\emat \\
\ns&\ns=\ns&\ns V^{\sss \top}\,A_{\scriptscriptstyle 0}^{\sss \top}\,X\,B\,(I_{\sss m}-
 R_{\sss X}^\dagger R_{\sss X}-R_{\sss X}^\dagger\,R)\,K_{\sss R}+\bmat{c} Q_{\sss 12} \\ Q_{\sss 22}
\emat \\
\ns&\ns=\ns&\ns V^{\sss \top}\,A_{\scriptscriptstyle 0}^{\sss \top}\,X\,B\,(I_{\sss m}-
 R_{\sss X}^\dagger R_{\sss X})\,K_{\sss R}+\bmat{c} Q_{\sss 12} \\ Q_{\sss 22}
\emat = V^{\sss \top}\,A_{\scriptscriptstyle 0}^{\sss \top}\,X\,B\,G_{\sss X}\,K_{\sss R}+\bmat{c} Q_{\sss 12} \\ Q_{\sss 22}
\emat.
\eeann
Recalling that 
 $\ima G_{\sss X}=\ker R_{\sss X}$, and that by virtue of (\ref{kercond0}) there holds $\ker R_{\sss X} \subseteq \ker (A_{\scriptscriptstyle 0}^{\sss \top}\,X\,B)$, we get $V^{\sss \top}\,A_{\scriptscriptstyle 0}^{\sss \top}\,X\,B\,G_{\sss X}\,K_{\sss R}=0$, from which {\bf (i)} immediately follows.
 To prove {\bf (ii)} we observe that 
 \bea
\label{psd}
\Pi_{\sss 1}=\bmat{ccc} I_{\sss n-\eta} & 0 & 0 \\ 0 & 0 & I_{\sss m} \emat
\bmat{c} A_{\sss V}^{\sss \top} \\ B_{\sss V}^{\sss \top} \emat
Q_{\sss V} [\begin{array}{cc} A_{\sss V} & B_{\sss V} \end{array}] \bmat{cc} I_{\sss n-\eta} & 0 \\  0 & 0 \\ 0 & I_{\sss m} \emat+\bmat{cc} 0 & 0 \\ 0 & R \emat \ge 0.
\eea
   In order to prove {\bf (iii)}, we first observe that in view of the previous considerations we have
 $X_{\sss V}= Q_{\sss V}+\bsmat \Delta_1 & 0 \\[1mm] 0 & 0 \esmat$. 
Substitution of this expression into (\ref{gdaretranf1}-\ref{gkercondtranf1}) yields
\beann
\bmat{cc}  \Delta_{\sss 1}  & 0   \\   0  & 0  \emat \ns&\ns = \ns&\ns 
\bmat{cc}   Q_{\sss 1}  & 0   \\   0  & 0  \emat+ \bmat{cc}   A_{\sss 1}^{\sss \top}\,\Delta_{\sss 1}\,A_{\sss 1}  & 0   \\   0  & 0  \emat \\
\ns&\ns \ns&\ns -\bmat{c} S_{\sss 1}+A_{\sss 1}^{\sss \top} \Delta_{\sss 1}\,B_{\sss 1} \\ 0 \emat (R_{\sss 1}+B_{\sss 1}^{\sss \top}\,\Delta_{\sss 1}\,B_{\sss 1})^{\dagger} [ \begin{array}{cc} S_{\sss 1}^{\sss \top}+B_{\sss 1}^{\sss \top} \,\Delta_{\sss 1}\,A_{\sss 1} & 0 \end{array} ],
\eeann
whose block in position (1,1) is exactly (\ref{redgdare2}).  We now prove that $\Delta_{\sss 1}$ satisfies
(\ref{redkercond2}).
Substitution of $X_{\scriptscriptstyle V}=Q_{\scriptscriptstyle V}+\bsmat \Delta_{\sss 1} & 0 \\[1mm] 0 & 0 \esmat$ into (\ref{gkercondtranf1}) gives
\beann
\ker (R_{\sss 1}+B_{\sss 1}^{\sss \top}\,\Delta_{\sss 1}\,B_{\sss 1}) \subseteq \ker \bmat{c} S_{\sss 1}+A_{\sss 1}^{\sss \top}\,\Delta_{\sss 1}\,B_{\sss 1}\\ \star \emat,
\eeann
from which (\ref{redkercond2}) immediately follows. \\
The second point can be proved by reversing these arguments along the same lines of the second part of the proof  of Theorem \ref{the1}.
\endproof

In view of {\bf (i)} of Theorem \ref{the2}, all solutions  of CGDARE($\Sigma$) coincide along 
$\gV \defi \ker \left( \bsmat I_{n-\eta} & 0 \\[1mm] 0 & 0 \esmat \,V^{\sss \top}\right)$. This means that given any two solutions $X$ and $Y$ of CGDARE($\Sigma$), we have $X|_{\gV}=Y|_{\gV}=Q_{\scriptscriptstyle 0}|_{\gV}$. \\[-3mm]

\begin{corollary}
\label{cor2}
The set $\gX$ of solutions of CGDARE($\Sigma$) is parameterized as the set of matrices
\[
X=V\,\bsmat \Delta_{\sss 1} && 0 \\[1mm] 0 && 0 \esmat \,V^{\sss \top}+Q_{\scriptscriptstyle 0}
\]
where $V=[\begin{array}{cc} V_{\sss 1} & V_{\sss 2} \end{array} ]$ is defined as in Theorem \ref{the2} and $\Delta_{\sss 1}$ is solution of (\ref{redgdare2}-\ref{redkercond2}).
\end{corollary}

{
\begin{remark}
{\em
In \cite{Ferrante-04} it is shown that if $X$ is a solution of 
DARE($\Sigma$) and we consider the associated solution $\Delta_{\sss 1}$ of the reduced DARE($\Sigma_{\sss 1}$), and if we denote by 
 $A_{\scriptscriptstyle X}$ and $A_{\sss \Delta_{\sss 1}}$  the associated closed-loop matrices, there holds
\bea
\label{hold}
V^{\sss \top} A_{\scriptscriptstyle X}\,V=\bmat{ll} A_{\sss \Delta_{\sss 1}} & 0 \\ \star & 0_{\sss \eta \times \eta} \emat.
\eea
This is a simple consequence of the fact that in the case of a solution $X$ of DARE($\Sigma$), the matrix $R_{\sss X}$ is invertible.
We now show via a simple example that this fact does not hold in general in the case of CGDARE($\Sigma$). Consider a Popov triple $\Sigma$ described by the matrices
\beann
A=\bsmat 
0 && 2 && 0 \\[1mm]
2 && 2 && 0 \\[1mm]
0 && 0 && -5 \esmat, \quad  
B=\bsmat 
-1 \\[1mm]
0 \\[1mm]
0  \esmat, \quad  
Q=\bsmat 
0 && 0 && 0 \\[1mm]
0 && 0 && 0 \\[1mm]
0 && 0 && 24 \esmat, \quad  
R=0, \quad  
S=\bsmat 
0  \\[1mm]
0   \\[1mm]
0 \esmat.
\eeann
In this case $A_{\sss 0}=A$ is invertible, and $A_{\sss 0}^{-1}\,B\,\ker R=\spanR \left\{ \bsmat 1 \\[1mm] -1 \\[1mm] 0\esmat \right\}$. Let $V_{\sss 2}=\bsmat -1/\sqrt{2} \\[1mm] 1/\sqrt{2} \\[1mm] 0 \esmat$ and 
$V=\bsmat -1/\sqrt{2} && 0 && -1/\sqrt{2} \\[1mm]
-1/\sqrt{2} && 0 && 1/\sqrt{2} \\[1mm]
0 && 1 && 0 \esmat$. Then, we compute
\beann
A_{\sss V} \ns&\ns = \ns&\ns  V^{\sss \top}\,A_{\sss 0} \,V=\bsmat 3 && 0 && -1 \\[1mm]
0 && -5 && 0 \\[1mm]
-1 && 0 && -1 \esmat, \quad 
B_{\sss V} = V^{\sss \top}\,B=\bsmat 1/\sqrt{2} \\[1mm] 0 \\[1mm]1/\sqrt{2} \esmat, \quad Q_{\sss V}=V^{\sss \top}\,Q_{\sss 0}\,V=\bsmat 0 && 0 && 0 \\[1mm]
0 && 24 && 0 \\[1mm]
0 && 0 && 0 \esmat,  \\
A_{\sss V}^{\sss \top}\,Q_{\sss V}\,A_{\sss V} \ns&\ns = \ns&\ns \bsmat 0 && 0 && 0 \\[1mm]
0 && 600 && 0 \\[1mm]
0 && 0 && 0 \esmat, \quad A_{\sss V}^{\sss \top}\,Q_{\sss V}\,B_{\sss V} =0,
\eeann
so that the matrices of the reduced CGDARE($\Sigma_{\sss 1})$ are
\beann
A_1 =\bsmat 3 && 0 \\[1mm] 0 && -5 \esmat, \quad B_{\sss 1}=\bsmat 1/\sqrt{2} \\[1mm] 0\esmat, \quad Q_1 =\bsmat 0 && 0 \\[1mm] 0 && 600 \esmat, \quad S_{\sss 1}=\bsmat 0 \\[1mm] 0\esmat,\quad R_{\sss 1}=0.
\eeann
A simple direct calculation shows that the only solution of this reduced CGDARE is $X_{\sss 1}=\bsmat 0 && 0 \\[1mm] 0 && -25 \esmat$. Thus, the only solution of the original CGDARE($\Sigma$) is $X=V\,\left(Q_{\sss V}+\bsmat X_1 && 0 \\[1mm] 0 && 0 \esmat \right)\,V^{\sss \top}=\bsmat 0 && 0 && 0 \\[1mm] 0 && 0 && 0 \\[1mm] 0 && 0 && -1 \esmat$. The corresponding closed-loop matrix coincides with $A$, i.e., $A_{\sss X}=A$. Now, 
\beann
V^{\sss \top}\,A_{\sss X}\,V=\bsmat 3 && 0 && -1 \\[1mm] 0 && -5 && 0 \\[1mm] -1 && 0 && -1 \esmat
\eeann 
This shows that neither of the two zero submatrices in the second block-column of (\ref{hold}) is zero in the general case of CGDARE($\Sigma$). While the submatrix in the upper left block of $A_{\sss X}$ still coincides with $A_{\sss \Delta_1}$, in the case of CGDARE($\Sigma$) it is also no longer true that 
the spectrum of $A_{\sss \Delta_1}$ is contained in that of $A_{\sss X}$. Indeed, in this case $\sigma(A_{\sss \Delta_1})=\{-5,3\}$ whereas 
$\sigma(A_{\sss X})=\{-5,1\pm\sqrt{5}\}$.
{This difference between DARE and CGDARE is related to the fact that in this generalized case the reduction can correspond simply to the singularity of $R$ which does not imply the singularity of $A_{\sss X}$ as discussed in Section \ref{sec2}. 
}}
\end{remark}
}

\begin{remark}
{\em
As for the reduction described in Theorem \ref{the1}, it may occur that, as a result of the reduction illustrated in Theorem \ref{the2}, $A_{\sss 1}-B_{\sss 1}\,R_{\sss 1}^{\dagger} S_{\sss 1}^{\sss \top}$ and/or $R_{\sss 1}$ be still singular. However, we have showed that $\Pi_{\sss 1}$ is symmetric and positive semidefinite. This means that if $A_{\sss 1}-B_{\sss 1}\,R_{\sss 1}^{\dagger} S_{\sss 1}^{\sss \top}$ is singular, we can repeat the reduction procedure described in Theorem \ref{the1}, while if $A_{\sss 1}-B_{\sss 1}\,R_{\sss 1}^{\dagger} S_{\sss 1}^{\sss \top}$ is non-singular but $R_{\sss 1}$ is singular, we can repeat the reduction procedure described in Theorem \ref{the2}. 
Since the order of the Riccati equation lowers at each reduction step, after at most $n$ steps, either we have computed the unique solution of the original CGDARE($\Sigma$), or we have obtained a symmetric Stein equation (which is linear), or we obtained a ``well-behaved'' DARE of maximally reduced order where the corresponding $R$ and $A-B\,R^{\dagger} S^{\sss \top}$ matrices are non-singular.
}
\end{remark}

\section{Numerical examples}

\begin{example}
\label{ex1}
{\em
Using the reduction techniques developed in the previous sections, we want to study the set of solutions of the CGDARE($\Sigma$) where $\Sigma$ is given by the matrices
\beann
A=\bsmat 
0 && -4 && 0 \\[1mm]
0 && 3 && 0 \\[1mm]
0 && 0 && -1 \esmat, \quad  
B=\bsmat 
0 && -1 \\[1mm]
3 && 0 \\[1mm]
0 && 0  \esmat, \quad  
Q=\bsmat 
1 && 0 && 0 \\[1mm]
0 && 0 && 0 \\[1mm]
0 && 0 && 0 \esmat, \quad  
R=\bsmat 
0&& 0 \\[1mm]
0 && 0 \esmat, \quad  
S=\bsmat 
0 && 0 \\[1mm]
0 && 0  \\[1mm]
0 && 0 \esmat.
\eeann
First notice that since $S$ is the zero matrix, $A_{\sss 0}$ and $Q_{\sss 0}$ coincide with $A$ and $Q$, respectively. Thus, in this case both $A_{\sss 0}$ and $R$ are singular. We begin with a reduction that corresponds to the singularity of $A_{\sss 0}$. Since $\ker A_{\sss 0}= \operatorname{span} \left\{ \bsmat 1 \\[1mm] 0 \\[1mm] 0 \esmat \right\}$, we can consider a basis matrix $U=[\begin{array}{c|c} U_{\sss 1} & U_{\sss 2} \end{array} ]$ given by
$U=\bsmat
0   &&  0 &&    1\\[1mm]
    -1  &&   0   &&  0\\[1mm]
         0  &&   1&&     0\esmat$,
so that 
\beann
A_{\sss U} = \bsmat
3   &&  0 &&   0\\[1mm]
 0  &&   -1   &&  0\\[1mm]
         4  &&   0&&     0\esmat, \quad \tilde{A} = \bsmat
3   &&  0 \\[1mm]
 0  &&   -1 \\[1mm]
         4  &&   0\esmat, 
 \quad B_{\sss U} = \bsmat
-3   &&  0 \\[1mm]
 0  &&    0\\[1mm]
         0  &&  -1\esmat,
          \quad Q_{\sss U} = \bsmat
0   &&  0 &&   0\\[1mm]
 0  &&   0   &&  0\\[1mm]
         0  &&   0&&     1\esmat.
         \eeann
         Thus, 
         \beann
         A_{\sss 1} =  \bsmat
3   &&  0 \\[1mm]
 0  &&   -1  \esmat, \quad 
          B_{\sss 1} =  \bsmat
-3   &&  0 \\[1mm]
 0  &&   0\esmat,
 \quad 
          S_{\sss 1} =  \bsmat
0  &&  -4 \\[1mm]
 0  &&   0\esmat,
 \quad  Q_{\sss 1} =  \bsmat
16 &&  0 \\[1mm]
 0  &&  0  \esmat,
  \quad  R_{\sss 1} =  \bsmat
0 &&  0 \\[1mm]
 0  &&  -1  \esmat.
 \eeann
 In view of Corollary \ref{cor1}, $X$ is a solution of CGDARE($\Sigma$) if and only if it can be written as
 \[
 X=Q_{\sss 0}+U\,\bsmat \Delta_{\sss 1} && 0 \\[1mm] 0 && 0 \esmat\,U^{\sss \top},
 \]
 where $\Delta_{\sss 1}$ is an arbitrary solution of (\ref{redgdare}-\ref{redkercond}). To maintain the notations as consistent as possible to those employed in 
 Section \ref{sec2b}, we define $\overline{A} \defi A_{\sss 1}$,  $\overline{B} \defi B_{\sss 1}$,  $\overline{Q} \defi Q_{\sss 1}$,  $\overline{S} \defi S_{\sss 1}$,  $\overline{R} \defi R_{\sss 1}$ and $\overline{X} \defi \Delta_{\sss 1}$.
 With this notation, (\ref{redgdare}-\ref{redkercond}) can be re-written as
 \bea
\overline{X}  \ns&\ns =\ns&\ns   \overline{A}_{\sss 0}^{\sss \top} \, \overline{X}  \,  \overline{A}_{\sss 0}  -   \overline{A}_{\sss 0}^{\sss \top}  \, \overline{X} \,  \overline{B} \,  (\overline{R}  +  \overline{B}^{\sss \top} \,  \overline{X}  \, \overline{B})^\dagger \overline{B}^{\sss \top} \,  \overline{X}\,\overline{A}_{\sss 0} +  \overline{Q}_{\sss 0} \label{es1a} \\
 &   &  \qquad \ker (\overline{R}+\overline{B}^{\sss \top}\,\overline{X}\,\overline{B}) \subseteq \ker (\overline{A}_{\sss 0}^{\sss \top}\,\overline{X}\,\overline{B}), \label{es1b}
 \eea
 where $\overline{A}_{\sss 0} = \overline{A}-\overline{B}\,\overline{R}^\dagger \,\overline{S}^{\sss \top}=\overline{A}$  and  $\overline{Q}_{\sss 0} = \overline{Q}- \overline{S}\, \overline{R}^\dagger \, \overline{S}^{\sss \top}= \bsmat
0 &&  0 \\[1mm]
 0  &&  0  \esmat$.
 Matrix $\overline{A}_{\sss 0}$ is invertible, whereas $\overline{R}$ is singular. Thus, we can apply the reduction procedure in Section \ref{sec2b} (we will employ the same notation used in Section \ref{sec2b}, with the only exception that all the letters will have a bar, to distinguish this second reduction from the first one). A simple calculation shows that $\ima (\overline{A}_{\sss 0}^{-1}\,\overline{B}\,\ker \overline{R})=\operatorname{span} \left\{ \bsmat 1 \\[1mm] 0 \esmat \right\}$. Thus, we can consider a basis matrix $V=[\begin{array}{c|c} V_{\sss 1} & V_{\sss 2} \end{array} ]$ given by
$V=\bsmat
0   &&      1\\[1mm]
    1&&     0\esmat$.
         Hence, we define $\overline{X}_{\sss V} \defi V^{\sss \top}\,\overline{X}\,V$ along with
         \beann
         \overline{A}_{\sss V}=V^{\sss \top}\,\overline{A}_{\sss 0}\,V=\bsmat -1 && 0 \\[1mm] 0 && 3 \esmat, \quad
          \overline{B}_{\sss V}=V^{\sss \top}\,\overline{B}=\bsmat 0 && 0 \\[1mm] -3 && 0 \esmat, \quad
            \overline{Q}_{\sss V}=V^{\sss \top}\,\overline{Q}_{\sss 0}\,V=\bsmat 0 && 0 \\[1mm] 0 && 0 \esmat, 
            \eeann
            so that 
      $\overline{A}_{\sss 1}  =-1$, $\overline{B}_{\sss 1}=\bsmat 0 && 0  \esmat$, $\overline{S}_{\sss 1}=\bsmat 0 && 0  \esmat$, $\overline{Q}_{\sss 1}  =0$, $\overline{R}_{\sss 1}=\bsmat 0 && 0 \\[1mm] 0 && 1 \esmat$.
           In view of Corollary \ref{cor2}, $\overline{X}$ is a solution of (\ref{es1a}-\ref{es1b}) 
            if and only if 
           \[
           \overline{X}=\overline{Q}_{\sss 0}+V\,\bsmat \overline{\Delta}_{\sss 1} && 0 \\[1mm] 0 && 0 \esmat\,V^{\sss \top}
           \]
           with $\overline{\Delta}_{\sss 1}$ being an arbitrary solution of 
            \bea
\overline{\Delta}_{\sss 1}  \ns&\ns =\ns&\ns   \overline{A}_{\sss 1}^{\sss \top} \overline{\Delta}_{\sss 1}  \overline{A}_{\sss 1}  -   \overline{A}_{\sss 1}^{\sss \top} \overline{\Delta}_{\sss 1} \overline{B}_{\sss 1} (\overline{R}_{\sss 1}  +  \overline{B}_{\sss 1}^{\sss \top} \overline{\Delta}_{\sss 1} \overline{B}_{\sss 1})^\dagger \overline{B}_{\sss 1}^{\sss \top} \overline{\Delta}_{\sss 1}\,\overline{A}_{\sss 1} +  \overline{Q}_{\sss 1} \label{es2a} \\
 &   &  \qquad \ker (\overline{R}_{\sss 1}+\overline{B}_{\sss 1}^{\sss \top}\,\overline{\Delta}_{\sss 1}\,\overline{B}_{\sss 1}) \subseteq \ker (\overline{A}_{\sss 1}^{\sss \top}\,\overline{\Delta}_{\sss 1}\,\overline{B}_{\sss 1}). \label{es2b}
 \eea
 We still have $\overline{R}_{\sss 1}$ singular, and $\overline{A}_{\sss 1}-\overline{B}_{\sss 1} \overline{R}_{\sss 1}^\dagger \overline{S}_{\sss 1}^{\sss \top}=\overline{A}_{\sss 1}$ is invertible. On the other hand, $\overline{B}_{\sss 1}\,\ker \overline{R}_{\sss 1}=\{0\}$, so that the reduction associated to the singularity of  $\overline{R}_{\sss 1}$ cannot be carried out. Using a change of coordinates in the input space given by $\Omega=\bsmat 0 && 1 \\[1mm] 1 && 0 \esmat$, we obtain
 \[
 \hat{R}_{\sss 1}=\Omega^{-1}\,\overline{R}_{\sss 1}\,\Omega=\bsmat 1 && 0 \\[1mm] 0 && 0 \esmat, \quad 
  \hat{B}_{\sss 1}=\overline{B}_{\sss 1}\,\Omega=\bsmat 0 && 0 \esmat,
  \]
so that $\hat{R}_{\sss 1,0}=1$ and  $\hat{B}_{\sss 1,0}=0$. Thus, (\ref{es2a}-\ref{es2b}) can be written in this basis as
  \bea
\overline{\Delta}_{\sss 1}  \ns&\ns =\ns&\ns   \overline{A}_{\sss 1}^{\sss \top} \overline{\Delta}_{\sss 1}  \overline{A}_{\sss 1}  -   \overline{A}_{\sss 1}^{\sss \top} \overline{\Delta}_{\sss 1} \hat{B}_{\sss 1,0} (\hat{R}_{\sss 1,0}  +  \hat{B}_{\sss 1,0}^{\sss \top} \overline{\Delta}_{\sss 1} \hat{B}_{\sss 1,0})^\dagger \hat{B}_{\sss 1,0}^{\sss \top} \overline{\Delta}_{\sss 1}\,\overline{A}_{\sss 1} +  \overline{Q}_{\sss 1} \label{es3a} \\
 &   &  \qquad \ker (\hat{R}_{\sss 1,0}+\hat{B}_{\sss 1,0}^{\sss \top}\,\overline{\Delta}_{\sss 1}\,\hat{B}_{\sss 1,0}) \subseteq \ker \overline{A}_{\sss 1}^{\sss \top}\,\overline{\Delta}_{\sss 1}\,\hat{B}_{\sss 1,0}. \label{es3b}
 \eea
which reduce to the trivial equation $\overline{\Delta}_{\sss 1}=\overline{\Delta}_{\sss 1}$ subject to the trivial constraint $\ker \bsmat 0 && 0 \\[1mm] 0 && 1 \esmat  \subseteq \ker \bsmat 0 && 0 \\[1mm] 0 && 0 \esmat$.
Any $\xi \defi \overline{\Delta}_{\sss 1}\in \real$ satisfies this reduced Riccati equation. Thus, {the solutions of  (\ref{es1a}-\ref{es1b}) are given by
$\overline{X}=V\,\bsmat \xi && 0 \\[1mm] 0 && 0 \esmat\,V^{\sss \top}=\bsmat 0 && 0 \\[1mm] 0 && \xi \esmat$, $\xi \in \real$, so that -- recalling that 
$Q_{\sss 0}=Q=\bsmat 1 && 0 && 0 \\[1mm] 0 && 0 && 0 \\[1mm] 0 && 0 && 0 \esmat$ and
 $U=\bsmat
0   &&  0 &&    1\\[1mm]
    -1  &&   0   &&  0\\[1mm]
         0  &&   1&&     0\esmat$ --
the set of solutions of the original CGDARE($\Sigma$) is parametrized by
\[
X=Q_{\sss 0}+U\,\bsmat 
0 && 0 & \Big| & 0 \\[-4mm]
0 && \xi & \Big| & 0 \\[-1mm]
\hline \\[-2mm]
0 && 0 & \Big| & 0 \esmat\,U^{\sss \top}=\bsmat 1 && 0 && 0 \\[1mm]
0 && 0 && 0 \\[1mm]
0 && 0 && \xi \esmat,\ \xi \in \real.
\]
}
 }
\end{example}

\begin{example}
{\em
Using the reduction techniques developed here, we want to study the set of solutions of the CGDARE($\Sigma$) where $\Sigma$ is given by the matrices
\beann
A=\bsmat 
4 && 0 && 0 \\[1mm]
-3 && 0 && 0 \\[1mm]
0 && 0 && -3 \esmat, \quad  
B=\bsmat 
3 && -5 \\[1mm]
1 && 1 \\[1mm]
0 && 0  \esmat, \quad  
Q=\bsmat 
3 && 0 && 0 \\[1mm]
0 && 0 && 0 \\[1mm]
0 && 0 && 16 \esmat, \quad  
R=\bsmat 
0&& 0 \\[1mm]
0 && 0 \esmat, \quad  
S=\bsmat 
0 && 0 \\[1mm]
0 && 0  \\[1mm]
0 && 0 \esmat.
\eeann
Since $S$ is the zero matrix, $A_{\sss 0}=A$ and $Q_{\sss 0}=Q$. Both $A_{\sss 0}$ and $R$ are singular. We begin with a reduction that corresponds to the singularity of $A_{\sss 0}$. Since $\ker A_{\sss 0}= \operatorname{span} \left\{ \bsmat 0 \\[1mm] 1 \\[1mm] 0 \esmat \right\}$, we can consider a basis matrix $U=[\begin{array}{c|c} U_{\sss 1} & U_{\sss 2} \end{array} ]$ given by
$U=\bsmat
1   &&  0 &&    0\\[1mm]
    0  &&   0   &&  1\\[1mm]
         0  &&   1&&     0\esmat$, so that 
\beann
A_{\sss U} = \bsmat
4   &&  0 &&   0\\[1mm]
 0  &&   -3   &&  0\\[1mm]
         -3  &&   0&&     0\esmat, \quad \tilde{A} = \bsmat
4  &&  0 \\[1mm]
 0  &&   -3 \\[1mm]
         -3  &&   0\esmat, 
 \quad B_{\sss U} = \bsmat
3   &&  -5 \\[1mm]
 0  &&    0\\[1mm]
         1  &&  1\esmat,
          \quad Q_{\sss U} = \bsmat
3   &&  0 &&   0\\[1mm]
 0  &&  16   &&  0\\[1mm]
         0  &&   0&&     0\esmat.
         \eeann
         Hence 
         \beann
         A_{\sss 1} =  \bsmat
4   &&  0 \\[1mm]
 0  &&   -3  \esmat, \quad 
          B_{\sss 1} =  \bsmat
3   &&  -5 \\[1mm]
 0  &&   0\esmat,
 \quad 
          S_{\sss 1} =  \bsmat
36  &&  -60 \\[1mm]
 0  &&   0\esmat,
 \quad  Q_{\sss 1} =  \bsmat
48 &&  0 \\[1mm]
 0  &&  144  \esmat,
  \quad  R_{\sss 1} =  \bsmat
27 &&  -45 \\[1mm]
 -45  &&  75  \esmat.
 \eeann
 In view of Corollary \ref{cor1}, $X$ is a solution of CGDARE($\Sigma$) if and only if it can be written as
 \[
 X=Q_{\sss 0}+U\,\bsmat \Delta_{\sss 1} && 0 \\[1mm] 0 && 0 \esmat\,U^{\sss \top},
 \]
 where $\Delta_{\sss 1}$ is an arbitrary solution of 
(\ref{redgdare}-\ref{redkercond}). As in Example \ref{ex1}, to maintain the notations as consistent as possible to those employed in 
 Section \ref{sec2b}, we define $\overline{A} \defi A_{\sss 1}$,  $\overline{B} \defi B_{\sss 1}$,  $\overline{Q} \defi Q_{\sss 1}$,  $\overline{S} \defi S_{\sss 1}$,  $\overline{R} \defi R_{\sss 1}$ and $\overline{X} \defi \Delta_{\sss 1}$.
 With this notation, (\ref{redgdare}-\ref{redkercond}) can be re-written as in (\ref{es1a}-\ref{es1b}) where
 \beann
 \overline{A}_{\sss 0} \ns&\ns=\ns&\ns \overline{A}-\overline{B}\,\overline{R}^\dagger \,\overline{S}^{\sss \top}= \bsmat
0 &&  0 \\[1mm]
 0  &&  -3  \esmat \quad \text{and} \quad
  \overline{Q}_{\sss 0} =  \overline{Q}- \overline{S}\, \overline{R}^\dagger \, \overline{S}^{\sss \top}= \bsmat
0 &&  0 \\[1mm]
 0  &&  144  \esmat.
 \eeann
Both $\overline{A}_{\sss 0}$ and $\overline{R}$ are singular. We can reapply the reduction procedure in Section \ref{secA0} (we will employ the same notation used in Section \ref{secA0}, with the only exception that all the letters will have a tilde, to distinguish this reduction from the first one). Now $\ker \overline{A}_{\sss 0}=\operatorname{span} \left\{ \bsmat 1 \\[1mm] 0 \esmat \right\}$. Thus, we can consider a basis matrix $\overline{U}=[\begin{array}{c|c}\overline{U}_{\sss 1} & \overline{U}_{\sss 2} \end{array} ]$ given by
$\overline{U}=\bsmat
0   &&      1\\[1mm]
    1&&     0\esmat$.
         Hence, we define $\overline{X}_{\sss \overline{U}} \defi \overline{U}^{\sss \top}\,\overline{X}\,\overline{U}$ along with
     $\overline{A}_{\sss \overline{U}}=\overline{U}^{\sss \top}\,\overline{A}_{\sss 0}\,\overline{U}=\bsmat -3 && 0 \\[1mm] 0 && 0 \esmat$, $\overline{B}_{\sss \overline{U}}=\overline{U}^{\sss \top}\,\overline{B}=\bsmat 0 && 0 \\[1mm] 3 && -5 \esmat$, $\overline{Q}_{\sss \overline{U}}=\overline{U}^{\sss \top}\,\overline{Q}_{\sss 0}\,\overline{U}=\bsmat 144 && 0 \\[1mm] 0 && 0 \esmat$.
We have thus obtained the matrices of the reduced-order Riccati equation
            \beann
           \overline{A}_{\sss 1}  =-3, \quad \overline{B}_{\sss 1}=\bsmat 0 && 0  \esmat, \quad \overline{S}_{\sss 1}=\bsmat 0 && 0  \esmat, \quad \overline{Q}_{\sss 1}  =1296, \quad  \overline{R}_{\sss 1}=\bsmat 0 && 0 \\[1mm] 0 && 0 \esmat.
           \eeann
           In view of Corollary \ref{cor2}, $\overline{X}$ is a solution of (\ref{es1a}-\ref{es1b}) if and only if 
           \[
           \overline{X}=\overline{Q}_{\sss 0}+\overline{U}\,\bsmat \overline{\Delta}_{\sss 1} && 0 \\[1mm] 0 && 0 \esmat\,\overline{U}^{\sss \top}
           \]
           with $\overline{\Delta}_{\sss 1}$ being an arbitrary solution of 
          (\ref{es2a}-\ref{es2b}).
      We still have $\overline{R}_{\sss 1}$ singular, and $\overline{A}_{\sss 1}-\overline{B}_{\sss 1} \overline{R}_{\sss 1}^\dagger \overline{S}_{\sss 1}^{\sss \top}=\overline{A}_{\sss 1}$ is invertible. On the other hand, $\overline{B}_{\sss 1}\,\ker \overline{R}_{\sss 1}=\{0\}$, so that the reduction associated to the singularity of  $\overline{R}_{\sss 1}$ cannot be carried out. Since 
 $\overline{R}_{\sss 1}$ is the zero matrix, and so is $\overline{B}_{\sss 1}$,  (\ref{es3a}-\ref{es3b}) can be written as the symmetric Stein equation 
  \beann
\overline{\Delta}_{\sss 1}  \ns&\ns =\ns&\ns   \overline{A}_{\sss 1}^{\sss \top} \overline{\Delta}_{\sss 1}  \overline{A}_{\sss 1}  +  \overline{Q}_{\sss 1} 
 \eeann
 subject to the trivial constraint $\ker (0) \subseteq \ker (0)$. This equation therefore reduces to
\beann
\overline{\Delta}_{\sss 1}  \ns&\ns =\ns&\ns    9\,\overline{\Delta}_{\sss 1}+1296 \label{redgdarees4}
  \eeann
  which admits the solution $\overline{\Delta}_{\sss 1}=-162$.
 Thus, 
 the matrix 
$\overline{X}=\overline{U}\,\bsmat -162 && 0 \\[1mm] 0 && 0 \esmat\,\overline{U}^{\sss \top}+\overline{Q}_{\sss 0}=\bsmat 0 && 0 \\[1mm] 0 &&-18 \esmat$ 
satisfies (\ref{es1a}-\ref{es1b}), and, recalling that 
$Q_{\sss 0}=Q=\bsmat 3 && 0 && 0 \\[1mm] 0 && 0 && 0 \\[1mm] 0 && 0 && 16\esmat$ and 
$U=\bsmat
1   &&  0 &&    0\\[1mm]
    0  &&   0   &&  1\\[1mm]
         0  &&   1&&     0\esmat$, we find
\[
X=Q_{\sss 0}+U\,\bsmat 
0 && 0 & \Big| & 0 \\[-4mm]
0 && -18 & \Big| & 0 \\[-1mm]
\hline \\[-2mm]
0 && 0 & \Big| & 0 \esmat\,U^{\sss \top}=\bsmat 3 && 0 && 0 \\[1mm]
0 && 0 && 0 \\[1mm]
0 && 0 && -2 \esmat,
\]
which is the only solution of the original CGDARE($\Sigma$).
 }
\end{example}

\section*{Concluding remarks}
{We have shown how a general CGDARE($\Sigma$) may be reduced to a well-behaved DARE($\Sigma$) of smaller order featuring a non-singular closed-loop matrix.
This reduction may be performed through repeated steps each of which may be easily implemented via robust linear algebraic routines thus providing an effective tool to deal with generalized Riccati equations in practical situations.}   

\end{document}